\documentclass[12pt]{amsart}
\usepackage{amssymb}
\begin{document}
%%%%%%%%%%%%%%%%%%%% Text italic %%%%%%%%%%%%%%%%%%%%%%%%%%%%
\theoremstyle{plain}
\newtheorem{thm}{Theorem}[section]
\newtheorem{theorem}[thm]{Theorem}
\newtheorem{lemma}[thm]{Lemma}
\newtheorem{corollary}[thm]{Corollary}
\newtheorem{proposition}[thm]{Proposition}
%%%%%%%%%%%%%%%%%%%% Text roman %%%%%%%%%%%%%%%%%%%%%%%%%%%%%
\theoremstyle{definition}
\newtheorem{notations}[thm]{Notations}
\newtheorem{question}[thm]{Question}
\newtheorem{remark}[thm]{Remark}
\newtheorem{remarks}[thm]{Remarks}
\newtheorem{definition}[thm]{Definition}
\newtheorem{claim}[thm]{Claim}
\newtheorem{assumption}[thm]{Assumption}
\newtheorem{properties}[thm]{Properties}
\numberwithin{equation}{thm}
%%%%%%%% Diagram macros, etc. %%%%%%%%%%%%%%%%%%%%%%%%%%%%%%%
\catcode`\@=11
% General macros
\def\opn#1#2{\def#1{\mathop{\kern0pt\fam0#2}\nolimits}}
\def\bold#1{{\bf #1}}%
\def\underrightarrow{\mathpalette\underrightarrow@}
\def\underrightarrow@#1#2{\vtop{\ialign{$##$\cr
 \hfil#1#2\hfil\cr\noalign{\nointerlineskip}%
 #1{-}\mkern-6mu\cleaders\hbox{$#1\mkern-2mu{-}\mkern-2mu$}\hfill
 \mkern-6mu{\to}\cr}}}
\let\underarrow\underrightarrow
\def\underleftarrow{\mathpalette\underleftarrow@}
\def\underleftarrow@#1#2{\vtop{\ialign{$##$\cr
 \hfil#1#2\hfil\cr\noalign{\nointerlineskip}#1{\leftarrow}\mkern-6mu
 \cleaders\hbox{$#1\mkern-2mu{-}\mkern-2mu$}\hfill
 \mkern-6mu{-}\cr}}}
% Rectangular Commutative diagrams
\let\amp@rs@nd@\relax
\newdimen\ex@
\ex@.2326ex
\newdimen\bigaw@
\newdimen\minaw@
\minaw@16.08739\ex@
\newdimen\minCDaw@
\minCDaw@2.5pc
\newif\ifCD@
\def\minCDarrowwidth#1{\minCDaw@#1}
\newenvironment{CD}{\@CD}{\@endCD}
\def\@CD{\def\A##1A##2A{\llap{$\vcenter{\hbox
 {$\scriptstyle##1$}}$}\Big\uparrow\rlap{$\vcenter{\hbox{%
$\scriptstyle##2$}}$}&&}%
\def\V##1V##2V{\llap{$\vcenter{\hbox
 {$\scriptstyle##1$}}$}\Big\downarrow\rlap{$\vcenter{\hbox{%
$\scriptstyle##2$}}$}&&}%
\def\={&\hskip.5em\mathrel
 {\vbox{\hrule width\minCDaw@\vskip3\ex@\hrule width
 \minCDaw@}}\hskip.5em&}%
\def\verteq{\Big\Vert&&}%
\def\noarr{&&}%
\def\vspace##1{\noalign{\vskip##1\relax}}\relax\iffalse{%
\fi\let\amp@rs@nd@&\iffalse}\fi
 \CD@true\vcenter\bgroup\relax\iffalse{%
\fi\let\\=\cr\iffalse}\fi\tabskip\z@skip\baselineskip20\ex@
 \lineskip3\ex@\lineskiplimit3\ex@\halign\bgroup
 &\hfill$\m@th##$\hfill\cr}
\def\@endCD{\cr\egroup\egroup}
% Horizontal arrows with "sliding" length
\def\>#1>#2>{\amp@rs@nd@\setbox\z@\hbox{$\scriptstyle
 \;{#1}\;\;$}\setbox\@ne\hbox{$\scriptstyle\;{#2}\;\;$}\setbox\tw@
 \hbox{$#2$}\ifCD@
 \global\bigaw@\minCDaw@\else\global\bigaw@\minaw@\fi
 \ifdim\wd\z@>\bigaw@\global\bigaw@\wd\z@\fi
 \ifdim\wd\@ne>\bigaw@\global\bigaw@\wd\@ne\fi
 \ifCD@\hskip.5em\fi
 \ifdim\wd\tw@>\z@
 \mathrel{\mathop{\hbox to\bigaw@{\rightarrowfill}}\limits^{#1}_{#2}}\else
 \mathrel{\mathop{\hbox to\bigaw@{\rightarrowfill}}\limits^{#1}}\fi
 \ifCD@\hskip.5em\fi\amp@rs@nd@}
\def\<#1<#2<{\amp@rs@nd@\setbox\z@\hbox{$\scriptstyle
 \;\;{#1}\;$}\setbox\@ne\hbox{$\scriptstyle\;\;{#2}\;$}\setbox\tw@
 \hbox{$#2$}\ifCD@
 \global\bigaw@\minCDaw@\else\global\bigaw@\minaw@\fi
 \ifdim\wd\z@>\bigaw@\global\bigaw@\wd\z@\fi
 \ifdim\wd\@ne>\bigaw@\global\bigaw@\wd\@ne\fi
 \ifCD@\hskip.5em\fi
 \ifdim\wd\tw@>\z@
 \mathrel{\mathop{\hbox to\bigaw@{\leftarrowfill}}\limits^{#1}_{#2}}\else
 \mathrel{\mathop{\hbox to\bigaw@{\leftarrowfill}}\limits^{#1}}\fi
 \ifCD@\hskip.5em\fi\amp@rs@nd@}
% Rectangular commutative diagrams with diagonal arows
\newenvironment{CDS}{\@CDS}{\@endCDS}
\def\@CDS{\def\A##1A##2A{\llap{$\vcenter{\hbox
 {$\scriptstyle##1$}}$}\Big\uparrow\rlap{$\vcenter{\hbox{%
$\scriptstyle##2$}}$}&}%
\def\V##1V##2V{\llap{$\vcenter{\hbox
 {$\scriptstyle##1$}}$}\Big\downarrow\rlap{$\vcenter{\hbox{%
$\scriptstyle##2$}}$}&}%
\def\={&\hskip.5em\mathrel
 {\vbox{\hrule width\minCDaw@\vskip3\ex@\hrule width
 \minCDaw@}}\hskip.5em&}
\def\verteq{\Big\Vert&}
\def\novarr{&}
\def\noharr{&&}
\def\SE##1E##2E{\slantedarrow(0,18)(4,-3){##1}{##2}&}
\def\SW##1W##2W{\slantedarrow(24,18)(-4,-3){##1}{##2}&}
\def\NE##1E##2E{\slantedarrow(0,0)(4,3){##1}{##2}&}
\def\NW##1W##2W{\slantedarrow(24,0)(-4,3){##1}{##2}&}
\def\slantedarrow(##1)(##2)##3##4{%
\thinlines\unitlength1pt\lower 6.5pt\hbox{\begin{picture}(24,18)%
\put(##1){\vector(##2){24}}%
\put(0,8){$\scriptstyle##3$}%
\put(20,8){$\scriptstyle##4$}%
\end{picture}}}
\def\vspace##1{\noalign{\vskip##1\relax}}\relax\iffalse{%
\fi\let\amp@rs@nd@&\iffalse}\fi
 \CD@true\vcenter\bgroup\relax\iffalse{%
\fi\let\\=\cr\iffalse}\fi\tabskip\z@skip\baselineskip20\ex@
 \lineskip3\ex@\lineskiplimit3\ex@\halign\bgroup
 &\hfill$\m@th##$\hfill\cr}
\def\@endCDS{\cr\egroup\egroup}
% Triangular commutative diagrams
\newdimen\TriCDarrw@
\newif\ifTriV@
\newenvironment{TriCDV}{\@TriCDV}{\@endTriCD}
\newenvironment{TriCDA}{\@TriCDA}{\@endTriCD}
\def\@TriCDV{\TriV@true\def\TriCDpos@{6}\@TriCD}
\def\@TriCDA{\TriV@false\def\TriCDpos@{10}\@TriCD}
\def\@TriCD#1#2#3#4#5#6{%
\setbox0\hbox{$\ifTriV@#6\else#1\fi$}
\TriCDarrw@=\wd0 \advance\TriCDarrw@ 24pt
\advance\TriCDarrw@ -1em
\def\SE##1E##2E{\slantedarrow(0,18)(2,-3){##1}{##2}&}
\def\SW##1W##2W{\slantedarrow(12,18)(-2,-3){##1}{##2}&}
\def\NE##1E##2E{\slantedarrow(0,0)(2,3){##1}{##2}&}
\def\NW##1W##2W{\slantedarrow(12,0)(-2,3){##1}{##2}&}
\def\slantedarrow(##1)(##2)##3##4{\thinlines\unitlength1pt
\lower 6.5pt\hbox{\begin{picture}(12,18)%
\put(##1){\vector(##2){12}}%
\put(-4,\TriCDpos@){$\scriptstyle##3$}%
\put(12,\TriCDpos@){$\scriptstyle##4$}%
\end{picture}}}
\def\={\mathrel {\vbox{\hrule
   width\TriCDarrw@\vskip3\ex@\hrule width
   \TriCDarrw@}}}
\def\>##1>>{\setbox\z@\hbox{$\scriptstyle
 \;{##1}\;\;$}\global\bigaw@\TriCDarrw@
 \ifdim\wd\z@>\bigaw@\global\bigaw@\wd\z@\fi
 \hskip.5em
 \mathrel{\mathop{\hbox to \TriCDarrw@
{\rightarrowfill}}\limits^{##1}}
 \hskip.5em}
\def\<##1<<{\setbox\z@\hbox{$\scriptstyle
 \;{##1}\;\;$}\global\bigaw@\TriCDarrw@
 \ifdim\wd\z@>\bigaw@\global\bigaw@\wd\z@\fi
 \mathrel{\mathop{\hbox to\bigaw@{\leftarrowfill}}\limits^{##1}}
 }
 \CD@true\vcenter\bgroup\relax\iffalse{\fi\let\\=\cr\iffalse}\fi
 \tabskip\z@skip\baselineskip20\ex@
 \lineskip3\ex@\lineskiplimit3\ex@
 \ifTriV@
 \halign\bgroup
 &\hfill$\m@th##$\hfill\cr
#1&\multispan3\hfill$#2$\hfill&#3\\
&#4&#5\\
&&#6\cr\egroup%
\else
 \halign\bgroup
 &\hfill$\m@th##$\hfill\cr
&&#1\\%
&#2&#3\\
#4&\multispan3\hfill$#5$\hfill&#6\cr\egroup
\fi}
\def\@endTriCD{\egroup}
%%%%%%%%%%%%%%%  End of diagram macros.  %%%%%%%%%%%%%%%%%%%%%%%%%
% Skriptbuchstaben
\newcommand{\sA}{{\mathcal A}}
\newcommand{\sB}{{\mathcal B}}
\newcommand{\sC}{{\mathcal C}}
\newcommand{\sD}{{\mathcal D}}
\newcommand{\sE}{{\mathcal E}}
\newcommand{\sF}{{\mathcal F}}
\newcommand{\sG}{{\mathcal G}}
\newcommand{\sH}{{\mathcal H}}
\newcommand{\sI}{{\mathcal I}}
\newcommand{\sJ}{{\mathcal J}}
\newcommand{\sK}{{\mathcal K}}
\newcommand{\sL}{{\mathcal L}}
\newcommand{\sM}{{\mathcal M}}
\newcommand{\sN}{{\mathcal N}}
\newcommand{\sO}{{\mathcal O}}
\newcommand{\sP}{{\mathcal P}}
\newcommand{\sQ}{{\mathcal Q}}
\newcommand{\sR}{{\mathcal R}}
\newcommand{\sS}{{\mathcal S}}
\newcommand{\sT}{{\mathcal T}}
\newcommand{\sU}{{\mathcal U}}
\newcommand{\sV}{{\mathcal V}}
\newcommand{\sW}{{\mathcal W}}
\newcommand{\sX}{{\mathcal X}}
\newcommand{\sY}{{\mathcal Y}}
\newcommand{\sZ}{{\mathcal Z}}
% Sonderbuchstaben mit Doppellinie
\newcommand{\A}{{\mathbb A}}
\newcommand{\B}{{\mathbb B}}
\newcommand{\C}{{\mathbb C}}
\newcommand{\D}{{\mathbb D}}
\newcommand{\E}{{\mathbb E}}
\newcommand{\F}{{\mathbb F}}
\newcommand{\G}{{\mathbb G}}
\newcommand{\HH}{{\mathbb H}}
\newcommand{\I}{{\mathbb I}}
\newcommand{\J}{{\mathbb J}}
\newcommand{\M}{{\mathbb M}}
\newcommand{\N}{{\mathbb N}}
\renewcommand{\P}{{\mathbb P}}
\newcommand{\Q}{{\mathbb Q}}
\newcommand{\R}{{\mathbb R}}
\newcommand{\T}{{\mathbb T}}
\newcommand{\U}{{\mathbb U}}
\newcommand{\V}{{\mathbb V}}
\newcommand{\W}{{\mathbb W}}
\newcommand{\X}{{\mathbb X}}
\newcommand{\Y}{{\mathbb Y}}
\newcommand{\Z}{{\mathbb Z}}
%%%%%%%%%%%%%%%%%%%%%%%%%%%%%%%%%%%%%%%%%%%%%%%%%%%%%%%%%%%%%%
\title[Brody hyperbolicity of moduli spaces]{On the Brody
hyperbolicity of moduli spaces for canonically
polarized manifolds}
%\title[]{} %mit Kurztitel in []
\author{Eckart Viehweg}
\address{Universit\"at GH Essen, FB6 Mathematik, 45117 Essen, Germany}
\email{ viehweg@uni-essen.de}
\thanks{This work has been supported by the ``DFG-Forschergruppe
Arithmetik und Geometrie'' and the ``DFG-Schwerpunktprogramm
Globale Methoden in der Komplexen Geometrie''}
\author[Kang Zuo]{Kang Zuo${}^*$}
\address{The Chinese University of Hong Kong, Department of Mathematics,
Shatin, Hong Kong}
\email{kzuo@math.cuhk.edu.hk}
\thanks{${}^*$ Supported by a ``Heisenberg-Stipendium'', DFG,
and partially by a grant from the Research
Grants Council of the Hong Kong
Special Administrative Region, China
(Project No. CUHK 4239/01P)}
\maketitle
%%%%%%%%%%% Introduction %%%%%%%%%%%%%%%%%%%%%%%%%%%%%%%%%%%
Given a polynomial $h$ let $\sM_h$ be the moduli functor of
canonically polarized complex manifolds with Hilbert polynomial
$h$. By \cite{Vie} there exists a coarse quasi-projective moduli
scheme $M_h$ for $\sM_h$, but in general $M_h$ will not carry a
universal family. Except for curves, there are no natural level
structures known, which can be added to enforce the existence of
fine moduli schemes. However, C. S. Seshadri and J. Koll\'ar
constructed finite coverings $Z \to M_h$ which are induced by a
``universal family'' in $\sM_h (Z)$ (see \cite{Vie}, 9.25).
Moreover, if a general element in $\sM_h ({\rm Spec} (\C))$ has no
non-trivial automorphism, then there exists an open subscheme
$M^{0}_{h} \subset M_h$ which carries a universal family. It is
the aim of this article to show that both, the coverings $Z$ and
the open subscheme $M^{0}_{h}$ are Brody hyperbolic. More general
we will show that the moduli stack $\sM_h$ is Brody hyperbolic,
in the following sense.

\begin{theorem} \label{0.1} Assume that for some quasi-projective
variety $U$ there exists a family $f: V \to U \in \sM_h (U)$ for which the
induced morphism $\varphi : U \to M_h$ is quasi-finite over its
image. Then $U$ is Brody hyperbolic; i.e. there are no non-constant
holomorphic maps $\gamma : \C \to U$.
\end{theorem}

Assume that the variety $U$ in \ref{0.1} is an open subvariety
of a projective $r$-dimensional manifold $Y$ with $B=Y\setminus U$ a normal
crossing divisor. We conjecture that the quasi-finiteness of $\varphi$
implies that $\Omega_Y^1(\log B)$ is weakly positive over some
open dense subset of $U$ (see definition \ref{3.1}) and that
$\kappa(\Omega_Y^r(\log B))=r$. \cite{Zuo} gives an affirmative
answer if for all the fibres of $V\to U$ the local Torelli
theorem holds true, and theorem \ref{0.1} adds some more evidence.

An algebraic version of \ref{0.1}, saying that for abelian varieties
$A$ or for $A = \C^*$ all algebraic morphisms $\gamma : A \to U$ have
to be constant, has been shown by S. Kov\'acs in \cite{Kov3} and
\cite{Kov1} (see also \cite{Mig}).

The non-existence of abelian subvarieties of moduli stacks
presumably can also be deduced from the bounds
for the degree of curves in moduli spaces (\cite{B-V},
\cite{V-Z} and \cite{Kov2}) by following the arguments used to
prove theorem 2.1 in \cite{Dem}.

Our arguments do not imply that the variety $U$ in
\ref{0.1} is hyperbolic in the sense of Kobayashi, except of
course if $U$ is a compact manifold and hence the Brody hyperbolicity
equivalent to the Kobayashi hyperbolicity. We will not
speculate about possible diophantine properties
of moduli schemes which conjecturally are related to
hyperbolicity (see \cite{Lan}).

A question similar to \ref{0.1} can be asked for moduli of
polarized manifolds, i.e. for the moduli functor of
pairs $(f:V\to U, \sH)$ where $f$ is a smooth projective
morphism with $\omega_F$ semi-ample for all fibres $F$ of $f$,
and where $\sH$ is fibrewise ample with Hilbert polynomial $h$.
Hence $\sP_h(U)$ is the set of such pairs, up to isomorphisms
and up to fibrewise numerically equivalence for $\sH$. By
\cite{Vie}, section 7.6, there exists a coarse quasi-projective
moduli scheme $P_h$ for $\sP_h$.

In \cite{V-Z} we have shown that for $U$ an elliptic curve,
or for $U=\C^*$ there are no non-isotrivial smooth
families $V\to U$, with $\omega_{V/U}$ relative semi-ample.
Being optimistic one could ask:
\begin{question} \label{0.2} Does the existence of some $(f: V \to
U,\sH) \in \sP_h (U)$ for which the induced morphism $\varphi :
U \to P_h$ is quasi-finite over its
image, imply that $U$ is Brody hyperbolic?
\end{question}

The methods used in this paper give an affirmative answer to
\ref{0.2} only under the additional assumption that
for some $\nu >0$ and for all fibres $F$ of $f$ the
$\nu$-canonical map $F \to \P(H^0(F,\omega_F^\nu))$
is smooth over its image. Except if $\omega_F^\nu=\sO_F$,
this additional assumption is by far too much to ask for, and
we do not consider this case in our article.\\

An outline of the content of this paper and a guideline to the
proof of \ref{0.1} will be given at the end of the
first section.\\

Luen Fai Tam and Ngaiming Mok introduced us to some of the
analytic methods used in this paper, and Luen Fai Tam checked a
preliminary version of section 7. We are grateful to both of
them for their interest and help.

%%%%%%%%%%%%%%%%%%%%%%% Sect. 1  %%%%%%%%%%%%%%%%%%%%%%%%%%%%%%%%%%
\section{A reformulation}

Theorem \ref{0.1} follows immediately from the next Propositions. In fact, if
there is a holomorphic map $\gamma : \C \to U$, we can replace $U$ by
the Zariski closure of $\gamma (\C)$, and the Proposition tells us that
the Zariski closure must be a point and hence $\gamma$ constant.

\begin{proposition} \label{1.1} Assume that for some $f: V \to U \in
\sM_h (U)$ the induced map $\varphi: U \to M_h$ satisfies
$$
\dim U = \dim \overline{\varphi (U)} > 0.
$$
Then there exists no holomorphic map $\gamma : \C \to U$ with
Zariski dense image.
\end{proposition}

The proposition \ref{1.1} is formulated in such
a way that, given a proper birational morphism $U' \to U$,
the assumptions allow to replace $f:V\to U$ by the fibre product
$f':V'=V\times_UU' \to U'$. We will call such a pullback family $f'$
a smooth birational model for $f$.

By the next lemma the conclusion in \ref{1.1} is compatible with
replacing $f$ by any smooth birational model.

\begin{lemma} \label{1.2} Let $\tau : U' \to U$ be a projective birational
morphism between quasi-projective varieties. Then a holomorphic map
$\gamma : \C \to U$ with Zariski dense image lifts to a holomorphic map
$\gamma' : \C \to U'$.
\end{lemma}

\begin{proof}
Let $U_0 \subset U$ be an open set with $\tau |_{\tau^{-1}
(U_0)}$ an isomorphism. $\gamma (\C)$ meets $U_0$, hence
$\gamma'$ exists on the complement of a discrete subset $A
\subset \C$. Let $\Delta$ be a small disk in $\C$, centered at $a
\in A$. The projective morphism $\tau$ factors through $U' \to U
\times \P^M$ for some $M$ and the composite $pr_2 \circ \gamma'
|_{\Delta^*}: \Delta^* \to \P^M$ is given by meromorphic
functions. Obviously it extends to a holomorphic map on
$\Delta$, and the image of the induced map $\Delta \to U\times
\P^M$ lies in $U'$.
\end{proof}

Using \ref{1.2} we will assume in the sequel that the quasi-projective
variety $U$ in \ref{1.1} is non-singular.\\

For the proof of \ref{1.1} we first gather and
generalize some methods of algebraic nature, in
particular the weak semi-stable reduction theorem of D. Abramovich
and K. Karu (see \cite{A-K}) and the positivity results for
direct images of certain sheaves (see \cite{Kaw}, \cite{Kol}, \cite{Vie2} and
\cite{Vie4}). In section 4 both will be applied to certain
product families, and the main result \ref{4.1} of this section
is quite similar to the one obtained by D. Abramovich in
\cite{Abr}. It will allow to replace the family $f:V\to U$ by
a smooth birational model of the $r$-fold product $f^r:V^r \to
U$ and to assume the stronger positivity properties stated in
\ref{4.3} and \ref{4.4}. Whereas the results of section 2 hold true for
arbitrary smooth projective morphisms, those of section 3
and 4 use the semi-ampleness of $\omega_F$ for all fibres $F$ of $f$.

Starting with section 5 we assume that contrary to \ref{1.1} or \ref{4.4}
there exists a holomorphic map $\gamma:\C \to U$
with dense image. In order to use covering constructions,
as we did in \cite{V-Z} for $\dim(U)=1$, we will
choose a hyperplane $H$ on $V$ whose discriminant
locus over $U$ is in general position with respect to
$\gamma(\C)$. At this point the ampleness of $\omega_F$ will be
needed.

In section 6 we will use the cyclic covering, obtained by taking
a root out of $H$ to compare and to study certain Higgs bundles
and their pullback to $\C$. The main properties are gathered in
\ref{6.2}. Finally section 7 contains some curvature estimates,
which show that the existence of $\gamma$, encoded in
lemma \ref{6.2}, contradicts the
Ahlfors-Schwarz lemma. The content of this section is influenced
by the work of J.-P. Demailly \cite{Dem}, S.S.-Y. Lu and S.-T.
Yau \cite{L-Y}, S.S.-Y. Lu \cite{Lu} and
Y.-T. Siu \cite{Siu} on hyperbolicity.

%%%%%%%%%%%%%%%%%%%%%%% Sect. 2 %%%%%%%%%%%%%%%%%%%%%%%%%%%%%
\section{Mild reduction}

Let $f: X \to Y$ be a morphism between projective manifolds with
connected general fibre. D. Abramovich and K. Karu constructed in
\cite{A-K} a generically finite proper morphism $Y' \to Y$ and a
proper birational map $Z' \to (X \times_Y Y'){\tilde{ \ }}$ such
that the induced morphism $g': Z' \to Y'$ is weakly semi-stable.
Here $\tilde{ \ }$ denotes the main component, i.e. the component dominant
over $X$,  We will not recall the definition of weak semi-stability,
but just list the main properties needed later.

\begin{definition} \label{2.1}
A morphism $g' : Z' \to Y'$ between projective varieties is called
mild, if
\begin{enumerate}
\item[a)] $g'$ is flat, Gorenstein with reduced fibres.
\item[b)] $Y'$ is non-singular and $Z'$ normal with at most rational
singularities.
\item[c)] Given a dominant morphism $Y'_1 \to Y'$ where $Y'_1$ has at
most rational Gorenstein singularities, $Z' \times_{Y'} Y'_1$ is
normal with at most rational singularities.
\item[d)] Let $Y'_0$ be an open subvariety of $Y'$, with
${g'}^{-1}(Y'_0)\to Y'_0$ smooth.
Given a non-singular curve $C'$ and a morphism
$\pi:C' \to Y'$ whose image meets $Y'_0$, the fibred product
$Z'\times_{Y'} C'$ is normal, Gorenstein with at most rational
singularities.
\end{enumerate}
\end{definition}
In \cite{A-K} the definition of a mild morphism just uses
the first three conditions, and by \cite{A-K}, 6.1 and
6.2, those hold true for weakly semi-stable morphisms.
As pointed out by K. Karu in \cite{Kar}, proof of 2.12, the
proof of the property c) carries over ``word by word'' to show
d). Hence d) holds true for weakly semi-stable morphisms as
well.

Hence starting with $f:X\to Y$, over some $Y'$, generically
finite over $Y$, one can find a
mild model of the pullback family, i.e. a mild morphism $g':Z'
\to Y'$ birational to $X\times_YY'\to Y'$. However it
might happen that one has to blow up the general fibre, and the
smooth locus of $g'$ will not be the pullback of the smooth
locus of $f$. Nevertheless, the existence of $g'$ will have
strong consequences for direct images of powers of dualizing
sheaves.

\begin{lemma} \label{2.2} Let $g' : Z' \to Y'$ be mild.
\begin{enumerate}
\item[i)] If $Y'' \to Y'$ is a dominant morphism between
manifolds, then $$pr_2 : Z' \times_{Y'} Y'' \to Y'' \mbox{ \ \
is mild.}
$$
\item[ii)] Let $g'' : Z'' \to Y'$ be a second mild morphism. Then
$$
(g' , g'') : Z ' \times_{Y'} Z'' \>>> Y'
\mbox{ \ \ is mild.}
$$
\item[iii)] For all $\nu \geq 1$ the sheaf $g'_* \omega^{\nu}_{Z'/Y'}$ is
reflexive.
\end{enumerate}
\end{lemma}

\begin{proof}
i) The property a) in \ref{1.2} is compatible with base change and
in c) one enforces the compatibility of b) with base change, as well.\\
ii) Since $Z''$ has rational Gorenstein singularities, the
property c) for $Z'$ implies that $Z' \times_{Y'} Z''$ has at
most rational Gorenstein singularities. The other properties
asked for in a) and b) are obvious. For c), remark that
$Z'' \times_{Y'} Y'_1$ is normal with rational Gorenstein
singularities and hence
$$
(Z'' \times_{Y'} Y'_1) \times_{Y'} Z' = (Z'' \times_{Y'} Z')
\times_{Y'} Y'_1
$$
has the same property. The same argument with $Y'_0$ replaced by
$C'$ gives d).

The sheaf $g'_* \omega^{\nu}_{Z'/Y'}$ is torsion free,
hence locally free outside of a closed codimension two
subvariety $T$ of $Y'$. Since $Z'$ is normal and equidimensional
over $Y'$, for $U_0 \subset Y'$ open and for $V_0={g'}^{-1}(U_0)$ one has
$$
H^0 (V_0, \omega^{\nu}_{Z'/Y'} ) = H^0 (V_0 \setminus {g'}^{-1}(T),
\omega^{\nu}_{Z'/Y'} ),
$$
and thereby
$$
H^0 (U_0, g'_* \omega^{\nu}_{Z'/Y'} ) = H^0 (U_0 \setminus T, g'_*
\omega^{\nu}_{Z'/Y'} ).
$$
So $g'_* \omega^{\nu}_{Z'/Y'}$ coincides with the maximal
extension of $g'_* \omega^{\nu}_{Z'/Y'}|_{Y'\setminus T}$ to $Y'$.
\end{proof}

Let $V\to U$ be any smooth projective morphism between
quasi-projective manifolds. We choose for $Y$ and $X$
projective non-singular compactifications, with $Y \setminus U$
and $X\setminus V$ normal crossing divisors, in such a way that
$V \to U$ extends to a morphism $f: X \to Y$.
If $g: Z' \to Y'$ denotes the weak semi-stable reduction, we
choose a birational morphism $\epsilon: Y_1 \to Y$ such that the
main component
$Y'_1=(Y' \times_{Y}Y_1)\tilde{ \ }$
is finite over $Y_1$. Let
$\Delta (Y'_1 / Y_1)$ denote the discriminant locus in $Y_1$ of
$Y'_1 \to Y_1$, and let $B_1 = Y_1 \setminus  \epsilon^{-1} (U)$ be the
boundary divisor. Blowing up a bit more we can assume that
$Y_1$ is non-singular, and that $\Delta (Y'_1 / Y_1) + B_1$ is a normal
crossing divisor.

By Y. Kawamata's covering construction (see \cite{Vie}, 2.6)
there exists a non-singular projective manifold $Y'_2$, finite
over $Y'_1$. In particular, there is a morphism $Y'_2 \to Y'$,
and by \ref{2.2}, a) the pullback of $Z'\to Y'$ is again mild.

Let us choose a desingularization $\psi: X_1 \to X \times_Y
Y_1$, such that
$$
(pr_2\circ \psi)^{*} (B_1+\Delta(Y'_1/Y_1))
$$
is a normal crossing divisor.

Changing the smooth birational model we may replace $U$ by its
pre-image in $Y_1$ and by abuse of notations rename $pr_2\circ
\psi : X_1 \to Y_1$ as $f: X \to Y$. We will also write
$Y'$ instead of $Y'_2$ and $Z'$ instead of $Z \times_{Y'} Y'_2$.
Doing so we reached the following situation:

\begin{lemma} \label{2.3} Any smooth projective morphism with
connected fibres has a smooth birational model $V\to U$ which
fits into a diagram of morphisms of normal varieties
\begin{equation}\label{2.3d}
\begin{CD}
V \> \subset >> X \< \tau' << X' \< \sigma << Z \< \rho << X'' \> \delta >>
Z' \\
\V VV \V V f V \V V f' V \V V g V \V V f'' V \V V g' V \\
U \> \subset >> Y \< \tau << Y' \< = << Y' \< = << Y' \> = >> Y'
\end{CD}
\end{equation}
with:
\begin{enumerate}
\item[i)] $Y$, $Y'$, $X$, $Z$ and $X''$ are non-singular projective
varieties.
\item[ii)] $\tau$ is finite and $X'$ is the normalization of $X
\times_Y Y'$.
\item[iii)] $\rho$ and $\delta$ are birational, and $\sigma$ is a
blowing up with center in the singular locus of $X'$.
\item[iv)] For $B = Y \setminus  U$ the divisors $B + \Delta (Y' / Y)$
and $f^* (B + \Delta (Y' / Y))$ are normal crossing divisors.
\item[v)] $g' : Z' \to Y'$ is mild.
\end{enumerate}
\end{lemma}

\begin{corollary} \label{2.4}
The conditions i) - v) stated in \ref{2.3} imply:
\begin{enumerate}
\item[vi)] $X'$ has rational singularities.
\item[vii)] For all $\nu \geq 1$ there exist isomorphisms
$$
g'_* \omega^{\nu}_{Z'/Y'} \> \simeq >> f''_*
\omega^{\nu}_{X''/Y'}\< \simeq << g_*\omega^{\nu}_{Z/Y'}.
$$
In particular, $g_*\omega^{\nu}_{Z/Y'}$ is a reflexive sheaf.
\item[viii)] For all $\nu \geq 1$ there exists an inclusion
$$
\iota : g_* \omega^{\nu}_{Z/Y'} \>>> \tau^* f_*
\omega^{\nu}_{X/Y},
$$
which is an isomorphism over $U$.
\item[ix)] For all $\nu \geq 1$ there exists some $N_\nu$ and an
invertible sheaf $\lambda_{\nu}$ on $Y$ with
$$
\tau^* \lambda_{\nu} \simeq \det (g_*
\omega^{\nu}_{Z/Y'})^{N_{\nu}}.
$$
\end{enumerate}
\end{corollary}

In part ix) the determinant of $g_* \omega^{\nu}_{Z/Y'}$ is
$i_* \det (g_* \omega^{\nu}_{Z/Y'} |_{Y \setminus  T})$ where $T$ is any
codimension two subvariety with $g_*
\omega^{\nu}_{Z/Y'}|_{Y\setminus T}$ locally free and $i : Y
\setminus  T \to Y$ the inclusion.

\begin{proof}
Since $\Delta (X'/X) \subset f^* \Delta (Y'/Y)$ are both normal
crossing divisors one obtains vi).

$Z'$ is normal with rational Gorenstein
singularities, hence $\delta^* \omega_{Z'/Y'} \subset
\omega_{X''/Y'}$ and $\omega^\nu_{Z'/Y'} =
\delta_*\delta^*\omega^\nu_{Z'/Y'} \subset \delta_* \omega^\nu_{X''/Y'}.$
The sheaf on the left hand side is invertible, and the one
on the right hand side torsionfree, and both coincide outside of
a codimension two subvariety. Hence they are equal and one
obtains the first isomorphism in vii). For the second one, one
repeats the argument for $\rho$ instead of $\delta$.
By \ref{2.2}, iii), all the three sheaves in vii) are reflexive.
Part viii) has been shown in \cite{Vie2}, 3.2 (see also \cite{Mor}, 4.10).

Let $B_{\nu}$ denote the zero divisor of $\det (\iota)$, hence
$$
\det (g_* \omega^{\nu}_{Z/Y'} ) \otimes \sO_{Y'} (B_{\nu}) =
\tau^* \det (f_* \omega^{\nu}_{X/Y}).
$$
In order to show that $B_{\nu}$ is the pullback of a
$\Q$-divisor on $Y$, we have to show, that the multiplicities
of two components of $B_\nu$ coincide, whenever both have the
same image in $Y$. To this aim, given any component
$\tilde{B}$ of $Y\setminus U$ consider a general curve $C$, which intersects
$\tilde{B}$ in some point $q$. Replacing $C$ by a neighborhood
of $q$ we will assume that this is the only intersection point.

Let us write $T_C=T\times_YC$, where $T$
stands for any of the varieties in the diagram (\ref{2.3d}).
Similarly, if $h:T\to T'$ is any of the morphism in the diagram
(\ref{2.3d}), $h_C$ will denote the restriction of $h$ to $T_C$.

By \ref{2.1}, d), the variety $Z'_C$ is again normal, Gorenstein
with at most rational singularities, and for $C$ sufficiently general
$X_C$ and $X''_C$ will be non-singular. Applying part viii) with
$Y$ replaced by $C$, one obtains a natural inclusion
\begin{equation}\label{incl2}
\iota_C: {g'_C}_* \omega^{\nu}_{{Z}'_C/{Y}'_C} \>>>
\tau_C^* ({f_C}_* \omega^{\nu}_{X_C/C}),
\end{equation}
and the zero divisor of $\det(\iota_C)$ is the restriction of
$B_\nu$ to $Y'_C$. In order to show ix), we just have to verify
that $B_\nu$ is the pullback of a $\Q$-divisor on $C$.

By \cite{KKMS} there exists a finite morphism $C' \to
C$, totally ramified in $q$, such that $X_C \times_{C} C'$
has a semi-stable model $S \to C'$.

By \ref{2.1}, d), the pullback of $Z'_C$ to some non-singular
covering of $C$ remains normal with rational Gorenstein
singularities. By flat base change (\ref{incl2}) is
compatible with further pullbacks.
Hence we may as well assume for a moment that ${Y}'_C \to C$
factors through $C'$. Then
$$
pr_1 : S' = S \times_{C'} Y'_C \>>> Y'_C \mbox{ \ \ and \ \ }
g'_C:Z'_C \>>> Y'_C
$$
are two flat Gorenstein morphism, $S'$ and $Z'_C$ are birational,
and both are normal with at most rational singularities.
Therefore, repeating the argument used to prove vii), one obtains
$${g'_C}_* \omega^{\nu}_{Z'_C/Y'_C} = pr_{1*}
\omega^{\nu}_{S'/Y'_C},
$$
and the divisor $B_{\nu}|_{Y'_C}$ is the pullback of a divisor
$\Pi$ on $C'$. Since $C'\to C$ is totally ramified in $q$,
the divisor $\Pi$ is itself the pullback of a $\Q$-divisor on
$C$.
\end{proof}
%%%%%%%%%%%%%%%%%%%%%%%% Sect.3 %%%%%%%%%%%%%%%%%%%%%%%%%
\section{Positivity of direct image sheaves}

As in \cite{Vie2} and \cite{Vie4} we use the following
convention: If $\sF$ is a coherent sheaf on a quasi-projective
normal variety $Y$, we consider the largest open subscheme $i: Y_1 \to
Y$ with $i^* \sF$ locally free. For
$$
\Phi = S^{\mu}, \ \ \
\Phi =\bigotimes^\mu \mbox{ \ \ \ or \ \ \ }\Phi = \det
$$
we define
$$
\Phi (\sF) = i_* \Phi (i^* \sF).
$$
\begin{definition} \label{3.1}
Let $\sF$ be a torsion free coherent sheaf on a quasi-projective
normal variety $Y$ and let $\sH$ be an ample invertible sheaf. Let $U \subset Y$ be an open subvariety.
\begin{enumerate}
\item[a)] $\sF$ is globally generated over $U$ if the natural
morphism
$$
H^0 (Y, \sF) \otimes \sO_Y \>>> \sF
$$
is surjective over $U$.
\item[b)] $\sF$ is weakly positive over $U$ if the restriction of $\sF$ to $U$ is locally free and if for all $\alpha >
0$ there exists some $\beta > 0$ such that
$$
S^{\alpha \cdot \beta} (\sF) \otimes \sH^{\beta}
$$
is globally generated over $U$.
\item[c)] $\sF$ is ample with respect to $U$ if there exists
some $\mu > 0$ such that
$$
S^{\mu} (\sF) \otimes \sH^{-1}
$$
is weakly positive over $U$.
\end{enumerate}
\end{definition}

The basic properties of weakly positive sheaves are listed in
\cite{Vie}, section 2.3. In particular, the definition of ``weak
positivity over $U$'' does not depend on the ample sheaf $\sH$
(\cite{Vie}, 2.14) and, if $\sF$ is weakly positive over $U$ and
$\sF \to \sG$ surjective over $U$, then $\sG$ is weakly positive
over $U$ (\cite{Vie}, 2.16). Moreover, weak positivity is a
local property. If for each point $u \in U$ there is a
neighborhood $U_0$ with $\sF$ weakly positive over $U_0$, then
$\sF$ is weakly positive over $U$.

By definition, most of the properties of weakly positive sheaves
$\sF$ carry over to sheaves which are ample over $U$.

\begin{lemma} \label{3.2} Let $\sH$ be an ample invertible sheaf
on $Y$. Then
$\sF$ is ample with respect to $U$, if its restriction to $U$ is locally free and if and only if for some
$\eta > 0$ there exists a morphism
$$
\bigoplus \sH \>>> S^{\eta} (\sF),
$$
surjective over $U$.
\end{lemma}

\begin{proof}
If $\sF$ is ample with respect to $U$, for all $\beta$
sufficiently large and divisible
$$
S^{2 \cdot \beta} (S^{\mu} (\sF)) \otimes \sH^{-2 \cdot \beta}
\otimes \sH^{\beta}
$$
is globally generated over $U$, as well as its
quotient sheaf
$$
S^{2 \cdot \beta \cdot \mu} (\sF) \otimes \sH^{-\beta}.
$$
We may assume that $\sH^{\beta-1}$ is very ample,
and we obtain the morphism asked for in \ref{3.2}.
On the other hand, if there is a morphism
$$
\bigoplus \sO_Y \>>> S^{\eta} (\sF) \otimes \sH^{-1},
$$
surjective over $U$, the sheaf $S^{\eta} (\sF) \otimes \sH^{-1}$
as a quotient of a weakly positive sheaf is weakly positive
over $U$.
\end{proof}

The basic methods to study positivity properties of direct
images are contained in \cite{Vie2}, \cite{Vie4}, \cite{Vie3} and
\cite{Vie}. Unfortunately in \cite{Vie2} and \cite{Vie4} we used
``weak positivity'' without specifying the open set, whereas in
\cite{Vie} we mainly work with smooth families or families
without non-normal fibres. So we have to recall some definitions
in this section, and we have to make the arguments carefully
enough to keep track of the open set $U$.

Let $f: X \to Y$ be a surjective projective morphism of
quasi-projective manifolds. We want to extend the constructions
from \cite{V-Z}, section 2, to the case $\dim(Y)>1$.

For an effective $\Q$-divisor $D\in {\rm Div}(X)$ the integral
part $[D]$ is the largest divisor with $[D] \leq D$.
For an effective divisor $\Gamma $ on $X$, and for
$N\in \N-\{0\}$ the algebraic multiplier sheaf is
$$
\omega_{X/Y} \left\{ \frac{-\Gamma}{N} \right\} = \tau_* \left(
\omega_{T/Y} \left( - \left[ \frac{\Gamma'}{N} \right] \right)
\right)
$$
where $\tau : T \to X$ is any blowing up with $\Gamma' = \tau_*
\Gamma$ a normal crossing divisor (see for example \cite{E-V},
7.4, or \cite{Vie}, section 5.3).

Let $F$ be a non-singular fibre of $f$.
Using the definition given above for $F$, instead of $X$, and for a divisor
$\Pi$ on $F$, one defines
$$
e (\Pi) = {\rm Min} \left\{ N \in \N \setminus  \{ 0 \} ; \ \omega_F
\left\{ \frac{- \Pi}{N} \right\} = \omega_F \right\} .
$$
By \cite{E-V} or \cite{Vie}, section 5.4, $e(\Gamma |_{F})$ is
upper semi-continuous, and there exists a neighborhood $V$ of $F$
with $e (\Gamma |_{V} ) \leq e (\Gamma |_{F})$. If $\sL$ is an
invertible sheaf on $F$, with $H^0(F,\sL)\neq 0$,
$$
e (\sL) = {\rm Max} \left\{e(\Pi) ; \ \Pi \mbox{ an effective
divisor and } \sO_F(\Pi)=\sL \right\} .
$$
\begin{proposition} \label{3.3}
Let $U \subset Y$ be an open subscheme, let $\sL$ be an
invertible sheaf, let $\Gamma$ be a divisor on $X$, and let $\sF$ be
a coherent sheaf on $Y$. Assume that, for some $N >0$, the
following conditions hold true:
\begin{enumerate}
\item[a)] $V = f^{-1} (U) \to U$ is smooth with connected fibres.
\item[b)] $\sF$ is weakly positive over $U$ (in particular $\sF
|_{U}$ is locally free).
\item[c)] There exists a morphism $f^* \sF \to \sL^N  ( -
\Gamma)$, surjective over $V$.
\item[d)] None of the fibres $F$ of $f: V \to U$ is contained in
$\Gamma$, and for all of them
$$
e (\Gamma |_{F}) \leq N.
$$
\end{enumerate}
Then $f_* (\sL \otimes \omega_{X/Y})$ is weakly positive over
$U$.
\end{proposition}

\begin{proof}
By \cite{Vie}, 5.23, the restriction of the sheaf $\sE = f_* (\sL \otimes
\omega_{X/Y})$ to $U$ is locally free. The verification of the
weak positivity will be done in several steps. Let us first show:
\begin{claim} \label{3.4}
In order to prove \ref{3.3} we are allowed to assume that $\sF$
is ample with respect to $U$.
\end{claim}
\begin{proof}
Let $\sH$ be a very ample sheaf on $Y$ and let $\rho: Y \to
\P^M$ be an embedding. For a general choice of the coordinate
planes $H_0, \ldots, H_M$, the intersection $H_i \cap (Y \setminus  U)$
is of codimension two in $Y$. We choose a codimension two
subscheme $T$ with $T \supset H_i \cap (Y\setminus U)$, for $i =0, \ldots
, M$. By definition, in order to show that $f_* (\sL \otimes
\omega_{X/Y})$ is weakly positive over $U$, we may replace $Y$
by $Y\setminus T$ and assume that $H_i \cap (Y\setminus U) =
\emptyset$. Moreover, for $T$ large enough $f$ will be flat. By the local
nature of weak positivity, it is sufficient to show that $f_*
(\sL \otimes \omega_{X/Y})$ is weakly positive over
$$
U_0 = U \setminus  \bigcup^{M}_{i=0} H_i.
$$
In fact, one can cover $U$ by such open sets,
for different choices of the coordinate planes.

Given $\alpha >0$ we choose $d = 1 + 2 \cdot \alpha$, and
consider the $d$-th power map
$$
\theta : \P^M \>>> \P^M \ \ \mbox{with} \ \ \theta (x_0, \ldots ,
x_M) = (x^{d}_{0}, \ldots , x^{d}_{M}).
$$
Let $Y'$ be the normalization of $\theta^{-1} (Y)$, and let
$\tau : Y' \to Y$ be the induced map. For the pullback $\sH'$ of
$\sO_{\P^1} (1)$ to $Y'$ one obtains $\tau^* \sH = {\sH'}^{d}.$

Leaving out codimension two subschemes in $Y$, not
meeting $U_0$, we may assume that $Y'$ is non-singular.
Then $X'= X \times_Y Y'$ is non-singular. In fact,
$f' : X' \to Y'$ is smooth over $\tau^{-1} (U)$ and $\tau' : X'
\to X$ is smooth over
$$
X \setminus  \bigcup^{M}_{i=0} f^{-1} (H_i).
$$
Let us choose $\sF' = \tau^* \sF \otimes {\sH'}^{N}$ and $\sL' = {\tau'}^{*}
\sL \otimes {f'}^{*} \sH'$. The sheaf $\sF'$ is ample with
respect to $U'_0=\tau^{-1}(U_0)$. Applying \ref{3.3} to $\sF'$ instead
of $\sF$ one finds
$$
f'_* (\sL' \otimes \omega_{X'/Y'}) = f'_* ({\tau'}^{*} \sL \otimes
\omega_{X'/Y'} ) \otimes \sH'
$$
to be weakly positive over $U'_0$. By flat base change, this sheaf
is isomorphic to
$$
\tau^* f_* (\sL \otimes \omega_{X/Y} ) \otimes \sH' = {\tau}^*(\sE)
\otimes \sH'.
$$
Hence for all $\beta$ sufficiently large and divisible, the sheaf
$$
S^{(2 \cdot \alpha) \cdot \beta} (\tau^* (\sE) \otimes \sH')
\otimes {\sH'}^{\beta} = \tau^*(S^{2 \cdot \alpha \cdot \beta}
(\sE)) \otimes {\sH'}^{(2\cdot \alpha+1)\cdot \beta}=\tau^*(S^{2
\cdot \alpha \cdot \beta}
(\sE) \otimes \sH^{\beta})
$$
is globally generated over $U'_0$. We obtain a morphisms
$$
\bigoplus \sO_{Y'} \>>> \tau^* (S^{2 \cdot \alpha \cdot
\beta} (\sE) \otimes \sH^{\beta}),
$$
surjective over $U'_0$, and
$$
\bigoplus \tau_* \sO_{Y'} \>>> S^{2 \cdot \alpha \cdot \beta} (\sE)
\otimes \sH^{\beta},
$$
surjective over $U_0$. For $\beta$ large enough, $\tau_* \sO_{Y'}
\otimes \sH^{\beta}$ is generated by global sections and hence
$S^{\alpha \cdot (2 \cdot \beta)} (\sE) \otimes \sH^{2 \cdot
\beta}$ is globally generated over $U_0$.
\end{proof}

\ref{3.4} allows to assume that $\sF$ is ample with respect to
$U$. Then the sheaf $\sL^{N \cdot \eta} (- \eta \cdot \Gamma)$
will be globally generated over $V$, for some $\eta \gg 0$.
Replacing $N$ by $N \cdot \eta$ and $\Gamma$ by $\eta \cdot
\Gamma$ we may as well assume that $\sL^N (- \Gamma)$ itself has
this property. From now on, this assumption will replace the
conditions b) and c) in \ref{3.3}.

Leaving out a codimension two subset of $Y$, not meeting $U$,
we will continue to assume that $f$ is flat. Let us fix some
non-singular compactification $\bar{Y}$ of $Y$ and a very ample
invertible sheaf $\bar{\sA}$ on $\bar{Y}$ such that
$\bar{\sA}^{\dim Y +1} \otimes \omega_{\bar{Y}}$ is ample.
We write $\sA=\bar{\sA}|_Y$ and $\sH = \sA^{\dim Y +1} \otimes
\omega_Y$.

\begin{claim} \label{3.5}
$\sE \otimes \sA^{\dim Y +1} \otimes \omega_Y$ is globally
generated over $U$.
\end{claim}

\begin{proof}
Let us choose a compactification $\bar{X}$ of $X$ such that $f$
extends to a morphism $\bar{f} : \bar{X} \to \bar{Y}$. Moreover
we choose $\bar{\sL}$ and $\bar{\Gamma}$ such that $\bar{\sL}^N
(- \bar{\Gamma})$ is again globally generated over $V$. Let
$\tau : X' \to \bar{X}$ be a blowing up, such that $\tau^*
\bar{\Gamma} = \Gamma'$ is a normal crossing divisor and let $f'
= \bar{f} \circ \tau$. The assumption d) in \ref{3.3} implies
that
$$
\sE' = f'_* \left(\tau^* \bar{\sL} \otimes \sO_{X'} \left( - \left[
\frac{\Gamma'}{N} \right] \right) \otimes \omega_{X'/\bar{Y}} \right)
\>>> \bar{f}_* (\bar{\sL} \otimes \omega_{\bar{X}/\bar{Y}})
$$
is an isomorphism over $U$. Hence it is sufficient to show that
$$
\sE' \otimes \bar{\sA}^{\dim Y +1} \otimes \omega_{\bar{Y}}
$$
is globally generated over $U$. Blowing up a bit more, and
enlarging $\Gamma'$ by adding components supported in $X' \setminus
\tau^{-1} (V)$, we can as well assume that $\tau^*(\bar{\sL})^N
\otimes\sO_{X'}(-\Gamma')$ is globally generated over $X'$. Under this
assumption \ref{3.5} has been shown in \cite{Vie}, 2.37, 2).
\end{proof}

To finish the proof, we consider for any $\alpha > 0$ the
$\alpha$-fold product
$$
X^{\alpha} = X \times_Y \ldots \times_Y X \mbox{ \ \
($\alpha$-times)}
$$
and $f^{\alpha} : X^{\alpha} \to Y$. Let
$\sigma: X^{(\alpha)} \to X^{\alpha}$ be a desingularization,
$f^{(\alpha)} = f^{\alpha} \circ \sigma$,
$$
\sL^{(\alpha)} = \sigma^* (\bigotimes^{\alpha}_{i=1} pr^{*}_{i}
\sL) \mbox{ \ \ \ and \ \ \ }
\Gamma^{(\alpha)} = \sigma^* (\sum^{\alpha}_{i=1} pr^{*}_{i}
\Gamma).
$$
$f^{(\alpha)} : X^{(\alpha)} \to Y$, $\sL^{(\alpha)}$
satisfies again the assumption a) in
\ref{3.3}. Moreover we assumed $\sL^N (- \Gamma)$ to be globally
generated over $V$, hence $\sL^{(\alpha)^N} (-
\Gamma^{(\alpha)})$ is globally generated over $V^r=V\times_U
\cdots \times_U V$. The assumption d) holds true for
$\Gamma^{(\alpha)}$ by:
\begin{claim} $e(\Gamma^{(\alpha)}|_{F^r})=e(\Gamma|_F)$
\end{claim}
\begin{proof}
The proof, similar to the one of 5.21 in \cite{Vie}, is by
induction on $r$. Obviously $e(\Gamma^{(\alpha)})\geq
e=e(\Gamma)$. Let $C$ be the support of the cokernel of
the inclusion
$$
\omega_{F^r}\left\{\frac{-\Gamma^{(\alpha)}|_{F^r}}{e}\right\}
\>>> \omega_{F^r}.
$$
Applying \cite{Vie}, 5.19, to the $i$-th projection $pr_i:F^r\to
F$, one finds subschemes $C_i$ of $F$ with $C=pr_i^{-1}(C_i)$.
Since this holds true for $i=1,\ldots, r$, $C$ must be empty.
\end{proof}
By \ref{3.5} the sheaf $f^{(\alpha)}_{*}
(\sL^{(\alpha)} \otimes \omega_{X^{(\alpha)}/Y} ) \otimes \sH$
is globally generated over $U$. Hence \ref{3.3} follows from the
next claim.

\begin{claim} \label{3.6} There exists a morphism
$$
f^{(\alpha)}_{*} (\sL^{(\alpha)} \otimes
\omega_{X^{(\alpha)}/Y} ) \>>> S^{\alpha} (f_* (\sL \otimes
\omega_{X/Y} )),
$$
surjective over $U$.
\end{claim}
\noindent
{\it Proof.} \
The natural morphism $\sigma_* \omega_{X^{(\alpha)}} \to
\omega_{X^{\alpha}}$ induces a morphism
$$
f^{(\alpha)}_{*} (\sL^{(\alpha)} \otimes
\omega_{X^{(\alpha)}/Y}) \>>> f^{\alpha}_{*}
((\bigotimes^{\alpha}_{i=1} pr^{*}_{i} \sL) \otimes
\omega_{X^{\alpha}/Y}),
$$
which is an isomorphism over $U$.
By flat base change, the right hand side is nothing but
$$
\bigotimes^{\alpha} f_* (\sL \otimes \omega_{X/Y}).
$$
\end{proof}

\begin{corollary} \label{3.7}
Let $f : X \to Y$ be a projective surjective morphism between
quasi-projective manifolds with connected general fibre. Assume
that for some open subscheme $U \subset Y$
$$
V = f^{-1} (U) \>>> U
$$
is smooth and that $\omega_{F_u}$ is semi-ample for all fibres $F_u=f^{-1}(u)$
with $u\in U$. Then $f_* \omega^{\nu}_{X/Y}$ is weakly positive over
$U$.
\end{corollary}

\begin{proof}
Using \ref{3.3} (with $\Gamma |_{V} =0$), one can copy the
arguments presented in the proof of Corollary 2.45 in
\cite{Vie} to obtain \ref{3.7} as a corollary to \ref{3.3}.
We leave this as an exercise since \ref{3.7} has been shown
under less restrictive assumptions in \cite{Vie3}, 3.7,
using different (and more complicated) arguments.
\end{proof}

\begin{remark}\label{basechange}
By \cite{Lev} the assumption ``$\omega_{F_u}$ is semi-ample for
all fibres $F_u$ with $u\in U$'' is equivalent to the
$f$-semi-ampleness of $\omega_{V/U}$. Hence for all $\nu$
sufficiently large and divisible, the natural morphism
$$
f^*f_*\omega_{X/Y}^\nu \>>> \omega_{X/Y}^\nu
$$
is surjective over $V$. In particular \ref{3.7} implies that
$\omega_{X/Y}$ is weakly positive over $V$.
\end{remark}
Let us end this section by stating a stronger positivity result.
Although it holds true by \cite{Kol} for arbitrary families of
manifolds of general type, we will just formulate it for
families with a semi-ample canonical sheaf. Recall that in
\cite{Vie2}, for a projective surjective morphism $f:X\to Y$ with
connected general fibre, we defined ${\rm Var}(f)$ to be the
smallest integer $\eta$ for which there exists a finitely
generated subfield $K$ of $\overline{\C(Y)}$ of transcendence degree $\eta$
over $\C$, a variety $F'$ defined over $K$, and a
birational equivalence
$$
X\times_Y {\rm Spec}(\overline{\C(Y)}) \sim F'\times_{{\rm Spec}(K)}
{\rm Spec}(\overline{\C(Y)}).
$$
\begin{theorem} \label{3.8} Under the assumptions made in \ref{3.7},
for all $\nu$ sufficiently large and divisible
$$
\kappa (\det (f_* \omega^{\nu}_{X/Y})) = {\rm Var}(f).
$$
\end{theorem}

\begin{proof} This has been shown in \cite{Vie4} if the
general fibres of $f$ are of general type, and in
\cite{Kaw} in general (see also \cite{Kol} or \cite{Vie3}).
\end{proof}

\begin{remark}\label{3.9} Let $f:V\to U$ be the morphism
considered in \ref{1.1}. Since $\omega_F$ is ample on the fibres
of $f$ we can replace the variety $F'$ in the definition of
${\rm Var}(f)$ by its image under a multicanonical map, hence
assume that it is also canonically polarized. One obtains a
morphism $\varphi':{\rm Spec}(K) \to M_h$ and
$K$ must contain the function field of
$\overline{\varphi(U)}_{\rm red}$.
In particular the assumption
$\dim(\overline{\varphi(U)})=\dim(U)$
implies that ${\rm Var}(f)=\dim(U)$.
\end{remark}

%%%%%%%%%%%%%%%%%%%%% Sect. 4 %%%%%%%%%%%%%%%%%%%%%%%%%%
\section{Products of families of canonically polarized
manifolds}

Let again $f: X \to Y$ be a surjective projective morphism
between quasi-projective manifolds with connected fibres and let
$U \subset Y$ be a non-empty open subvariety such that
$$
f: V = f^{-1} (U) \>>> U
$$
is smooth, and such that $\omega_{V/U}$ is $f$-semi-ample.

In \cite{V-Z}, 2.7, we have shown that for curves $Y$, the
ampleness of $\det (f_* \omega^{\nu}_{X/Y})$ implies the
ampleness of $f_* \omega^{\nu}_{X/Y}$, for $\nu \geq 2$. In
\cite{Vie}, 6.22, one finds a similar statement over $U$. In
order to extend the latter to $Y$, one would like to control the
``non-local free locus'' of $f_* \omega^{\nu}_{X/Y}$. This could
be done by using natural compactifications of moduli spaces, but
those only exist for curves, for surfaces of general type, or by
\cite{Kar} under strong assumptions on the existence of minimal
models.

Fortunately, the mild reduction of D. Abramovich and K. Karu can
serve as a substitute, using in particular the reflexivity
of the sheaves in \ref{2.4} vii).

We will assume in the sequel, that $\dim(U)={\rm Var}(f)$, and that
$V \to U$ fits into the diagram considered in \ref{2.3}. Since
$Y'$ is finite over $Y$ one finds ${\rm Var}(g)={\rm
Var}(f)=\dim(Y')$, and \ref{3.8} implies that $\det (g_*
\omega^{\nu}_{Z/Y'})$ is big for all $\nu$ sufficiently large
and divisible. We choose
such $\nu \geq 3$, and we will assume in addition that
$$
f^{*} f_* \omega^{\nu}_{V/U} \>>> \omega^{\nu}_{V/U}
$$
and the multiplication morphisms
$$
S^\beta(f_*\omega^{\nu}_{V/U}) \>>> f_*\omega^{\beta\cdot\nu}_{V/U}
$$
are surjective, for all $\beta$. Define
$$
e = {\rm Max} \{ e (\omega^{\nu}_{F}) ; \ F \mbox{ a fibre of} \ V
\to U \} .
$$
By \ref{2.4}, ix), there is an invertible sheaf $\lambda_{\nu}$
on $Y$ and some $N_{\nu} \in \N$ with
$$
\tau^* \lambda_{\nu} = \det (g_* \omega^{\nu}_{Z/Y'}
)^{N_{\nu}} .
$$
Writing $B=Y\setminus U$ for the boundary divisor,
let us fix an ample invertible sheaf $\sA$, such that
$\sA(-B)$ is ample. Since
$$
\kappa(\lambda_{\nu})=\kappa(\det (g_*
\omega^{\nu}_{Z/Y'}))=\dim(Y),
$$
there exists some $\eta > 0$ and some effective divisor $D$, with
$\lambda^{\eta}_{\nu} = \sA (D)$. Replacing $N_{\nu}$ by
some multiple we can assume
$$
\det (g_* \omega^{\nu}_{Z/Y'})^{N_{\nu}} = \tau^* \sA
(D)^{\nu \cdot (\nu -1) \cdot e}.
$$
Define $r_0 = {\rm rank} (f_* \omega^{\nu}_{X/Y})$ and $r =
N_{\nu} \cdot r_0$.

\begin{proposition} \label{4.1}
Let $X^{(r)}$ denote a desingularization of the $r$-th fibre
product $X \times_Y \ldots \times_Y X$ and let $f^{(r)} : X^{(r)}
\to Y$ be the induced morphism. Then for all $\beta$
sufficiently large and divisible the sheaf
$$
f^{(r)}_*(\omega^{\beta\cdot\nu}_{X^{(r)}/Y}) \otimes
\sA^{-\beta\cdot\nu\cdot (\nu-2)} \otimes
\sO_{Y} (-\beta\cdot\nu\cdot (\nu -1) \cdot D)
$$
is globally generated over $U$ and
$$
\omega^{\beta\cdot\nu}_{X^{(r)}/Y} \otimes f^{(r)*}
(\sA^{-\beta\cdot\nu\cdot(\nu-2)} \otimes \sO_{Y} (-
\beta\cdot\nu\cdot(\nu -1) \cdot D))
$$
is globally generated over $V^r = f^{(r)-1} (U)$.
\end{proposition}
\begin{proof}
We use again the notations from \ref{2.3}. By \ref{2.2}, ii),
mildness of a morphism is compatible with fibre products, hence
$$
{g'}^r : {Z'}^r = Z' \times_{Y'} \ldots \times_{Y'} Z' \to Y'
$$
is again mild.

For the normalization ${X'}^{(r)}$ of $X^{(r)} \times_Y Y'$ we
choose a desingularization ${Z}^{(r)}$, with centers in the
singular locus of ${X'}^{(r)}$, and a non-singular blowing up
${X''}^{(r)}$ which dominates both, $Z^{(r)}$ and ${Z'}^r$. We
obtain again a diagram
$$
\begin{CD}
V^r \> \subset >> X^{(r)} \< \tau^{(r)} << {X'}^{(r)} \<
\sigma^{(r)} << {Z}^{(r)} \<
\rho^{(r)} << {X''}^{(r)} \> \delta^{(r)} >> {Z'}^r \\
\V VV \V V f^{(r)} V \V VV \V g^{(r)} VV \V V {f''}^{(r)} V \V V {g'}^r V \\
U \> \subset >> Y \< \tau << Y' \< = << Y' \< = << Y' \> = >> Y'
\end{CD}
$$
which satisfies the assumptions made in \ref{2.3}. One finds for
all integers $\mu \geq 0$
\begin{equation} \label{equa1}
{g'}^{r}_{*} \omega^{\mu}_{{Z'}^r / Y'} = \bigotimes^r g'_*
\omega^{\mu}_{Z'/Y'}.
\end{equation}
In fact, by flat base change and by the projection formula, both
sheaves coincide over the largest subvariety of $Y'$, where
$g'_*\omega^\mu_{Z'/Y'}$ is locally free. By definition, the
right hand side of (\ref{equa1}) is the reflexive hull of the
tensor product on this subscheme, and by \ref{2.2}, iii) the left
hand side is reflexive, hence both are equal. Corollary \ref{2.4}
implies:

\begin{claim} \label{4.2} \ \\
\begin{enumerate}
\item[a)] ${g}^{(r)}_{*} \omega^{\mu}_{{Z}^{(r)}/Y'}$ is
reflexive and there is an isomorphism
$$
{g}^{(r)}_{*} \omega^{\mu}_{{Z}^{(r)}/Y'} \simeq \bigotimes^r
g_* \omega^{\mu}_{Z/Y'}.
$$
\item[b)] There is an inclusion
$$
{g}^{(r)}_{*} \omega^{\mu}_{{Z}^{(r)}/Y'} \>>> \tau^*
f^{(r)}_{*} \omega^{\mu}_{X^{(r)}/Y}
$$
which is an isomorphism over $U'$.
\end{enumerate}
\end{claim}

\begin{proof}
b) and the first part of a) is nothing but \ref{2.4}, viii) and
vii). For the second part of a), \ref{2.4}, vii) allows to
replace the left hand side by ${g'}^{r}_{*}
\omega^{\mu}_{{Z'}^{r}/Y'}$, the right hand side by
$\bigotimes^r g'_* \omega^{\mu}_{Z'/Y'}$, and to apply
(\ref{equa1}).
\end{proof}

By construction $g^{(r)}:Z^{(r)} \to Y'$ is smooth over
$U' = \tau^{-1} (U)$ and
$$
{g^{(r)}}^{-1}(U')= {V'}^r = V^r\times_UU'.
$$
Now we play the usual game. For the integer $\nu \geq 3$ chosen
above, and for $r_0 = {\rm rank} (g_*\omega^{\nu}_{Z/Y'})$,
there is a natural inclusion
\begin{equation} \label{incl0}
\det (g_* \omega^{\nu}_{Z/Y'}) \>>> \bigotimes^{r_0} (g_*
\omega^{\nu}_{Z/Y'} )
\end{equation}
which locally splits over the open set $Y'_1$, where
$g_*\omega^{\nu}_{Z/Y'}$ is locally free, in particular over
$U'$. By the choice of $r$ one obtains an inclusion
\begin{equation} \label{incl1}
\tau^* \sA (D)^{\nu \cdot (\nu -1) \cdot e} \>>> \bigotimes^r (g_*
\omega^{\nu}_{Z/Y'}) = {g}^{(r)}_*\omega^{\nu}_{{Z}^{(r)}/Y'},
\end{equation}
again locally splitting over $U'$.
In fact, the splitting inclusions in (\ref{incl0}) and (\ref{incl1})
exist over $Y'_1$, and since the sheaves on the right hand sides
are reflexive they extend to $Y'$.

For $\omega = \omega_{{Z}^{(r)}/Y'}$ and $\sA' = {g}^{(r)*}
\tau^* \sA (D)^{\nu}$ consider $\sL = \omega \otimes {\sA'}^{-1}$.
By (\ref{incl1}) $\omega^{\nu} \otimes {\sA'}^{-(\nu -1) \cdot
e}$ has a section whose zero divisor $\Gamma$ does not contain
a whole fibre over $U'$, and
$$
\sL^{\nu \cdot (\nu -1) \cdot e} = \omega^{\nu \cdot (\nu -1)
\cdot e - \nu^2} \otimes \omega^{\nu^2} \otimes {\sA'}^{-\nu
\cdot (\nu -1) \cdot e} = \omega^{\nu \cdot (\nu -1) \cdot e -
\nu^2} \otimes \sO_{{Z}^{(r)}} (\nu \cdot \Gamma ).
$$
All fibres of $V^r \to U$ are of the form
$$
F^r=F\times \cdots \times F
$$ and \cite{Vie}, 5.21, implies that
$$
e (\Gamma |_{F^r}) \leq e (\omega^{\nu}_{F^r}) = e
(\omega^{\nu}_{F} ) \leq e.
$$
So $e( \nu \cdot \Gamma |_{F^r}) \leq \nu \cdot e$ and for $N =
\nu \cdot e$ the assumption b) in \ref{3.3} holds true. By corollary
\ref{3.7} the sheaf ${g}_* \omega^{\nu \cdot (\nu \cdot (e-1)-e)}$
is weakly positive over $U'$. Since
$$
{g}^{(r)*} {g}^{(r)}_{*} \omega^{\nu \cdot (\nu \cdot (e-1)
-e)} \>>> \omega^{\nu \cdot (\nu \cdot (e-1)
-e)} = \sL^{\nu \cdot (\nu -1) \cdot e} \otimes
\sO_{{Z}^{(r)}} (- \nu \Gamma).
$$
is surjective over $V^r$ we can apply \ref{3.3} (for
$\sL^{\nu-1}$ instead of $\sL$) and obtain the
weak positivity of
$$
{g}^{(r)}_{*} (\sL^{\nu -1} \otimes \omega_{{Z}^{(r)}/Y'} )
= {g}^{(r)}_{*} ( \omega^{\nu}_{{Z}^{(r)}/Y'}) \otimes
\tau^* \sA (D)^{-\nu \cdot (\nu -1)}
$$
over $U'$. Since $g^{(r)}_*\omega^{\beta\cdot\nu}_{Z^{(r)}/Y'}$
is reflexive, one has the multiplication morphism
$$
\mu_\beta:S^\beta(g^{(r)}_*\omega^\nu_{Z^{(r)}/Y'}) \>>>
g^{(r)}_*\omega^{\beta\cdot\nu}_{Z^{(r)}/Y'}.
$$
By \ref{4.2}, a), the left hand side is
$S^\beta(\bigotimes^rg_*\omega^\nu_{Z/Y'})$
whereas the right hand side is
$\bigotimes^rg_*\omega^{\beta\cdot\nu}_{Z/Y'}$,
hence the assumption on the surjectivity of the multiplication
morphism carries over, and $\mu_\beta$
is surjective over $U'$. Since
$$
{g}^{(r)}_{*} (\omega^{\nu}_{{Z}^{(r)}/Y'}) \otimes \tau^* \sA
(D)^{-\nu \cdot (\nu -1)}
$$
is weakly positive over $U'$, for all $\beta$
sufficiently large and divisible
$$
S^{\beta}({g}^{(r)}_{*} ( \omega^{\nu}_{{Z}^{(r)}/Y'}) \otimes
\tau^* \sA (D)^{-\nu \cdot (\nu -1)})\otimes \tau^*\sA^{\beta}
$$
as well as
$$
{g}^{(r)}_{*} ( \omega^{\beta\cdot\nu}_{{Z}^{(r)}/Y'}) \otimes
\tau^* \sA (D)^{-\beta\cdot\nu \cdot (\nu -1)}\otimes \tau^*\sA^{\beta}
$$
will be globally generated over $U'$. By \ref{4.2}, b), one has
a morphism
$$
\tau_*\sO_{Y'}\otimes\sA^{\beta\cdot(\nu-1)} \>>>
\tau_*\tau^*({f}^{(r)}_{*} \omega^{\beta\cdot\nu}_{{X}^{(r)}/Y}) \otimes
\sA (D)^{-\beta\cdot\nu \cdot (\nu -1)}\otimes \sA^{\beta\cdot\nu},
$$
surjective over $U$. Although the sheaf ${f}^{(r)}_{*}
\omega^{\beta\cdot\nu}_{{X}^{(r)}/Y}$ is not necessarily
reflexive, the finiteness of $\tau$ allows to apply the projection
formula, and to obtain thereby a morphism
$$
\tau_*\sO_{Y'}\otimes\sA^{\beta\cdot(\nu-1)}\>>>
{f}^{(r)}_{*}( \omega^{\beta\cdot\nu}_{{X}^{(r)}/Y}) \otimes
\sA^{-\beta\cdot\nu \cdot (\nu -2)}\otimes \sO_Y(-\beta\cdot\nu
\cdot (\nu -1) \cdot D),
$$
surjective over $U$. For $\beta$ large enough,
the sheaf on the left hand side will be generated
by global sections, hence for those $\beta$ the sheaf on the
right hand side is globally generated over $U$. Since we assumed
$$
f^{(r)*} f^{(r)}_{*} \omega^{\nu}_{X^{(r)/Y}} \>>>
\omega^{\nu}_{X^{(r)}/Y}
$$
to be surjective over $V$, the same holds true for $\nu$
replaced by $\beta\cdot\nu$, and
$$
\omega^{\beta\cdot\nu}_{X^{(r)}/Y} \otimes f^{(r)*} (\sA^{-
\beta\cdot\nu \cdot (\nu - 2)} \otimes \sO_{X^{(r)}} (- \beta\cdot\nu
\cdot (\nu -1) \cdot D))
$$
is globally generated over $V^r$.
\end{proof}
From now on, we will forget the original morphism $f$ and
work only with the morphism $f^{(r)}$. To keep notations
as simple as possible, we allow ourselves to change them
again. Doing so, we can restate the results of the sections
2, 3 and 4 in the following way:
\begin{corollary} \label{4.3}
Let $\tilde{U}$ be a quasi-projective manifold and let
$\tilde{f}:\tilde{V}\to \tilde{U}$ be a smooth projective surjective
morphism with connected fibres,
with ${\rm Var}(\tilde{f})=\dim(\tilde{U})$ and with $\omega_{\tilde{F}}$
semi-ample for all fibres $\tilde{F}$ of $\tilde{f}$.

Then there exist a proper birational morphism $U \to \tilde{U}$,
a projective compactification $Y$ of $U$, a projective morphism
$f:X\to Y$, an invertible sheaf $\sA$ and an
effective divisor $D$ on $Y$, such that for all $\nu$
sufficiently large and divisible one has:
\begin{enumerate}
\item[a)] $f:V=f^{-1}(U) \to U$ is smooth with connected fibres.
\item[b)] $X$ and $Y$ are projective manifolds, and $X\setminus V$ and
$B=Y\setminus U$ are normal crossing divisors.
\item[c)] $\sA$ is ample, and $D\geq B$.
\item[d)] $f_*(\omega_{X/Y}^\nu) \otimes \sA(D)^{-\nu}$ is
globally generated over $U$.
\item[e)] $\omega_{X/Y}^\nu \otimes f^* \sA(D)^{-\nu}$ is
globally generated over $V$.
\end{enumerate}
\end{corollary}
\begin{proof}
By \ref{2.3} we find a smooth birational model $f:V\to U$ of
$\tilde{f}:\tilde{V} \to \tilde{U}$ which
fits into the diagram in \ref{2.3}. We may
replace $V\to U$ by $V^r \to U$, and apply \ref{4.1}.
The property a) and b) obviously hold true. Since we assumed
$\sA(-B)$ to be ample and $\nu\geq 3$, for the invertible
$\sA'=\sA^{\nu-2}(-B)$ and for the divisor
$$
D'=(\nu-1)\cdot D + B
$$
one obtains property c) and by \ref{4.1} d) and e).
\end{proof}

If one starts with any smooth morphism in \ref{1.1}, one knows
by \ref{3.9} that the variation is maximal. \ref{1.2} allows to
blow up the base, hence \ref{4.3} allows to replace the original
morphism by a new one, satisfying the assumptions a) - e).
Thereby the proposition \ref{1.1} and hence theorem \ref{0.1} are
immediate consequences of the next proposition, which will be
shown at the end of section 7.
\begin{proposition}\label{4.4} Given $U$, let $f:X\to Y$ be a
projective surjective morphism, satisfying the conditions a) - e)
in \ref{4.3} for some $\nu$, $\sA$ and $D$. Assume moreover that
$n=\dim(F)$ is even, that $r=\dim(U) \geq 1$, and
that $\omega_F^\nu$ is very ample for all fibres $F$ of $V\to U$.
Then there exists no holomorphic map $\gamma:\C \to U$ with
dense image.
\end{proposition}
%%%%%%%%%%%%%%%%%%%%%%%% Sect. 5 %%%%%%%%%%%%%%%%%%%%%
\section{Construction of cyclic coverings}

Starting from a morphism $f:X\to Y$ satisfying the assumptions
in \ref{4.4} for an invertible sheaf $\sA$, a
divisor $D$ and a natural number $\nu$ let us consider
$$
\sL = \omega_{X/Y}\otimes f^* \sA(D)^{-1}.
$$
Blowing up $X$ with centers outside of $V$ we may assume that
the global sections of $\sL^\nu$ generate an invertible sheaf
$\sH$. If $E$ denotes the divisor on $X$ with $\sH(E)=\sL$,
then $E$ has support in $X\setminus V$, hence it is a normal
crossing divisor.

Let us assume there exists a holomorphic map $\gamma: \C \to U$
with dense image, contrary to \ref{4.4}. In this section we will
choose some divisor and some cyclic covering of $X$, depending
on $\gamma$ and finally this construction will help to show that
such a holomorphic map can not exist.

By \ref{4.3}, d), we have for some $\ell$ a morphism
$\displaystyle\bigoplus^{\ell+1} \sO_Y \to f_* \sL^\nu$,
surjective over $U$, and by \ref{4.3}, e) the induced morphisms
$$
\bigoplus^{\ell+1} \sO_X \>>> f^*f_* \sL^\nu \>>> \sL^\nu
$$
are both surjective over $V$. By assumption one obtains
embeddings
$$
V \>>> \P=\P(f_*\sL^\nu|_V)\>>> \P^\ell \times U.
$$
The projection to $\P^\ell$ extends to the morphism
$\pi:X \to \P^\ell$, defined by the sections of the sheaf
$\sH$. For all hyperplanes $H$ in $\P^\ell$ one
has $$\sL=\sO_X(\pi^*(H)+E).$$

Let $\check{\P}^\ell$ denote the dual projective space.
For a hyperplane $H \subset \P^\ell$, we will write $[H]\in
\check{\P}^\ell$ for the corresponding point. For each $u\in U$
and for $F_u=f^{-1}(u)$ the set of all $[H]\in\check{\P}^\ell$
with $F_u\cap H$ non-singular and not equal to $F_u$ is open.
Let $S_u$ denote the complement. By \cite{D-K}, XVII, 3.2, for general
points $u$ of $(\ell-1)$-dimensional components of $S_u$, the
intersection $F_u\cap H$ will have just one ordinary
double point of type $A_1$, i.e. a singularity given
locally analytic as the zero set of the equation
$x_1^2+\ldots+x_n^2$ in $\C^n$. Hence the locus $T_u$, consisting of hyperplanes
$H$ with $F_u\cap H$ having other types of singularities or with
$F_u \subset H$, is of codimension at least two in $\check{\P}^\ell$.

As in \cite{D-K}, XVII, 6.1 those properties can also be considered
in families, and the corresponding sets depend algebraically on
the parameter. In particular,
$$
S=\{([H],u) ; \ F_u \subset H \mbox{ \ or \ } F_u\cap H \mbox{ \
singular}\}
$$
is a closed subset of $\check{\P}^\ell \times U$. Let us choose
a codimension $2$ closed subscheme $T$ of $\check{\P}^\ell \times
U$, contained in $S$ such that $S\setminus T$ is non-singular,
of pure codimension one, and
$$
S\setminus T \subset \{([H],u) ; \ F_u \not\subset H \mbox{ \
and \ } F_u\cap H \mbox{ \ has one double point of type \ } A_1 \}.
$$
Given $[H]\in \check{\P}^\ell$ let $S_H$ and $T_H$ be the
intersection of $\{[H]\}\times
U$ with $S$ and $T$, respectively.
\begin{lemma}\label{5.1} There exists some
$[H]\in\check{\P}^\ell$ such that $T_H \cap \gamma(\C) =
\emptyset$, such that $S_H$ meets $\gamma(\C)$ transversally,
and such that $\pi^*(H)$ is non-singular and $\pi^*(H)+E$ a
normal crossing divisor.
\end{lemma}
Here ``$S_H$ meets $\gamma(\C)$ transversally'' means that for
a local section $\sigma$ of $\sO_U$ with zero set $(S_H)_{\rm red}$,
the holomorphic function $\gamma^*(\sigma)$ has zeros of order one.
\begin{proof} $\gamma:\C \to U$ induces a holomorphic map
$$
\tilde{\gamma}: \check{\P}^\ell \times \C \>>> \check{\P}^\ell\times U.
$$
Since $\tilde{\gamma}$ is holomorphic, $\Delta^{(1)}=\tilde{\gamma}^{-1}(T)$
is a complex subspace of $\check{\P}^\ell\times \C$.
Let $\Delta^{(2)}$ be the complex subspace of
$\tilde{\gamma}^{-1}(S)$ given locally by the following condition.
Let $\sigma$ be a local equation of $S$ on
${\check{\P}^\ell\times U}$. Then $\Delta^{(2)}$ is the analytic
subspace of the zero set of $\tilde{\gamma}^*\sigma$, where the
multiplicity of $\tilde{\gamma}^*\sigma$ is larger than or equal
to two. We choose $\Delta=\Delta^{(1)} \cup \Delta^{(2)}$.

By \cite{G-R}, page 172, $\Delta$ has a decomposition
$$
\Delta=\bigcup_{i\in I} \Delta_i
$$
in irreducible components. The index set $I$ is countable, since
each point $p\in \C$ has a small neighborhood $U(p)$ such that
$\check{\P}^\ell\times U(p)$ meets only finitely many of those
components. As usual,
$$
\dim(\Delta)={\rm Max}  \{ \dim(\Delta_i); \ i\in I\}.
$$
\begin{claim}\label{5.2} \ \ $\dim(\Delta) \leq \ell - 1$.
\end{claim}
\begin{proof} If $\gamma$ is not an embedding of a small
neighborhood of a point $p\in \C$, then
$$
\check{\P}^\ell\times \{p\}\cap \Delta^{(2)}
$$
consists of all hyperplanes $H$ passing through $p$, and its
dimension is $\ell -1$. The set of those points is discrete.
For all other points $p$ and for all
components $\Delta_i$ of $\Delta^{(2)}$ one has
$$
\dim(\check{\P}^\ell\times \{p\}\cap \Delta_i) \leq \ell - 2.
$$
In fact, let $U(p)$ denote a sufficiently small neighborhood of
$p$. A general $[H]\in\check{\P}^\ell$ does not pass through
$\gamma(p)$, and for those who do, the intersection is
transversal, except for all $[H]$ in a codimension $2$ subset
of $\check{\P}^\ell$.

If $\Delta_i$ is one of the components of $\Delta^{(1)}$, then
for all $p\in\C$
$$
\dim(\Delta_i\cap \check{\P}^\ell\times \{p\})=\dim(T\cap
\check{\P}^\ell\times \{ \gamma(p)\}) \leq \ell -2.
$$
In both cases, if $\Delta_i$ is a component of $\Delta$ with
$\Delta_i \subset \check{\P}^\ell\times\{p\}$,
we are done. Otherwise choose for $j=1, \ 2$ two points $p_j\in \C$ with
$$
\check{\P}^\ell\times\{p_j\}\cap\Delta_i\neq \emptyset.
$$
Then $\check{\P}^\ell\times\{p_1\}\cap\Delta_i$ is not dense in
$\Delta_i$. Obviously, the dimension of
$\check{\P}^\ell\times\{p_1\}\cap\Delta_i$ is larger than or
equal to $\dim(\Delta_i)-1$, and by Ritt's lemma (\cite{G-R},
page 102) both must be equal. Hence
$$
\dim(\Delta_i)=\dim(\check{\P}^\ell\times\{p_1\}\cap\Delta_i)+1
\leq \ell -1.
$$
\end{proof}
\begin{claim}\label{5.3} \ \ $pr_1(\Delta)$ does not contain an
open analytic subset $W\subset \check{\P}^\ell$.
\end{claim}
\begin{proof}
We will show \ref{5.3} by induction on $\ell$, just using
\ref{5.2} but not the definition of $\check{\P}^\ell$ as a dual
projective space.
If $\ell=1$, the set $pr_1(\Delta)$ is countable.

In general, if $W\subset pr_1(\Delta)$ we choose a point
$p\in\C$, such that none of the countably many components of
$\Delta$ is contained in $\check{\P}^\ell\times\{p\}$. Moreover,
for each $i\in I$, we choose a point $q_i\in pr_1(\Delta_i)$. Let
$H\simeq \check{\P}^{\ell-1}$ be a hyperplane, passing through
$p$ but not containing any of the points $q_i$. Then, for each
component $\Delta_i$, the intersection $\Delta_i\cap H\times \C$
can not be dense in $\Delta_i$ and
$$
\dim(\Delta_i\cap H\times \C) < l-1.
$$
Hence
$$
\dim(\Delta\cap H\times \C) \leq \ell-2=\dim(H) -1,
$$
and since $W\cap H$ is an open analytic subset of $H$, contained
in $pr_1(\Delta\cap H\times \C)$, this contradicts the
induction hypotheses.
\end{proof}
Recall that we assumed that the global sections of $\sL$
generate the invertible subsheaf $\sH$ of $\sL$. In particular,
$$
H^0(X,\sH)=H^0(X,\sL^\nu)=H^0(\P^\ell,\sO_{\P^\ell}(1)),
$$
and for $[H]$ in some Zariski open subscheme $\check{\P}^\ell_0$
of $\check{\P}^\ell$ the preimage $\pi^*(H)$ will be
non-singular and $\pi^*(H)+E$ a normal crossing divisor. By
\ref{5.3} we can find points $[H]$ in
$\check{\P}_0^\ell\setminus pr_1(\Delta)$, and for all of them the
properties asked for in \ref{5.1} hold true.
\end{proof}
From now on $H$ is fixed, and we write $T=B\cup T_H$ and $S$ for
the closure of $S_H$ in $Y$.
We will not use anymore the fact that $T_H$ is of codimension one,
and in the next step we will replace $Y$ by a blow up with
the centers partly contained in $T_H$.
\begin{lemma}\label{5.4} Assume that, contrary to \ref{4.4},
there exists $\gamma:\C\to U$ with a dense image. Then we may
assume in addition to \ref{4.3}, a), b), d) and e), that there exists a general
section of $\sL^\nu=\omega^\nu_{X/Y}\otimes f^*\sA(D)^{-\nu}$ with zero
divisor $H+E$, and divisors $S$ and $T$ in $Y$ such that:
\begin{enumerate}
\item[i)] $S\cap U$ is dense in $S$ and $S+T$ and $f^*(S+T)$ are
normal crossing divisors.
\item[ii)] $X\to Y$ and $H\to Y$ are both smooth over
$U_0=Y\setminus (S\cup T)$.
\item[iii)] The fibres of $H\to Y$ over $Y_0=Y\setminus T$ are
reduced with at most an ordinary double point.
\item[iv)] $\gamma(\C)\cap T = \emptyset$.
\item[v)] $H$ is non-singular, and $f(E)$ is contained in $B$.
\item[vi)] $\sA$ is semi-ample, ample with respect to $Y_0$, and $D\geq B$.
\end{enumerate}
\end{lemma}
\begin{proof} All this can be done by blowing up $Y$ in centers
not contained in $\gamma(\C)$ and replacing $f:X\to Y$ by a
desingularization of the pullback family.
\end{proof}
The section of $\sL^\nu$ with zero divisor $H+E$ gives rise to a
cyclic covering $\psi':Z'\to X$ (see for example \cite{E-V},
Section 3). The condition \ref{5.4}, ii), implies that
$$
g: Z_0={\psi'}^{-1}f^{-1}(U_0) \>>> U_0
$$
is smooth, hence it gives rise to a variation of Hodge
structures $\V_0=R^ng_{*}\C_{Z_0}$.
\begin{lemma}\label{5.5} The monodromy of
$\V_0=R^n g_{*}\C_{Z_0}$ around the components of $S$ is
finite.
\end{lemma}
\begin{proof} Here we will use the assumption, that the
dimension $n$ of the fibres of $f$ is even.

A general curve $C$ meets $S$ transversally. Replacing $C$
by some open subset, we can assume that for a given component
$S_i$ of $S$
$$
C\cap (S\cup T)=C\cap S_i= \{p\}.
$$
The restriction
$$
\psi_C: Z_C=Z'\times_YC \>>> X_C=X\times_YC
$$
of the finite morphism $\psi':Z'\to X$
is a cyclic covering of order $\nu$, totally ramified along
$H_C=H\times_YC$. By the definition of $S$ and $T$, the fibre
$H_p=H_C\cap F_p$ has one singular point $q$, and we can choose
locally analytic parameters $t$ in a neighborhood of $p\in C$ and $t,x_1, \ldots, x_n$
in a neighborhood of $q\in X_C$ such that $H_C$ is the zero-set
of $\sum_{i=1}^n x_i^2+t$ near $q$. Then locally near $\psi_C^{-1}(q)$
the covering $Z_C$ is given by the equation
$$
\sum_{i=1}^n x_i^2+t+z^{\nu}.
$$
So $g^{-1}(p)$ has one isolated singularity, a double point of type $A_{\nu-1}$.
As well known (see \cite{Loo}, p. 132, for example), in even dimension
the local monodromy group of an $A_{\nu-1}$ singularity is finite,
and as in \cite{D-K} or \cite{Loo}, p. 41, one obtains the same for the global monodromy.
\end{proof}

%%%%%%%%%%%%%%%%%%%%%%%% Sect. 6 %%%%%%%%%%%%%%%%%%%%
\section{Higgs bundles}

\begin{notations}\label{5.6}
In this section we assume that $f:X\to Y$ satisfies the
conditions stated in \ref{4.4}, except possibly that $\sA$ is not
ample, but only semi-ample and big. For the given holomorphic map
$\gamma:\C \to U$ we assume moreover the existence of the
divisors $S$, $T$, $H$ and $E$, satisfying the conditions in
\ref{5.4}.

We define $\Delta=f^*(T)$ and $\Sigma=f^*(S)$.
Recall, that the original boundary divisor $B$ is contained in
$T$. So the non-reduced components of $\Delta$ or the
components of $\Delta+\Sigma$, mapping to codimension two
subvarieties of $Y$, are all supported in $f^{-1}B$.

Let $\delta:X'\to X$ be a blowing up of $X$ with centers in
$\Delta+\Sigma$ such that $H'+\Delta'+\Sigma'$ is a
normal crossing divisor, where
$\Delta'=\delta^*\Delta$, $\Sigma'=\delta^*\Sigma$ and where
$H'$ is the proper transform of $H$.
For $\sL=\omega_{X/Y}\otimes f^*\sA(D)^{-1}$, we write $\sL'=\delta^*\sL$.
For $E'=\delta^*(H+E)-H'$ on finds ${\sL'}^\nu=\sO_{X'}(H'+E')$.

Let $g:Z_0 \to U_0$ be the fibre space, considered at the end of the
last section, obtained by restricting the cyclic covering
$\psi':Z'\to X$, given by the divisor $H+E$ in \ref{5.4}.
We choose $Z$ to be a desingularization of the normalization of
the fibre product $X'\times_X Z'$. Let us denote the induced
morphisms by
$$
\begin{CD}
Y\< g << Z \> \delta' >> Z'\\
\V = V V \V \psi V V \V V \psi' V \\
Y \< f' << X' \> \delta >> X.
\end{CD}
$$
Finally we write $\Pi=g^{-1}(S\cup T)$, and identify $Z_0$ with
$Z\setminus\Pi$.

In the sequel we will write $T_*(-\log \bullet)$ for the dual of
$\Omega^1_*(\log \bullet)$.
\end{notations}

By \cite{Del}, for all $k\geq 0$ the local constant system $R^kg_*\C_{Z_0}$
gives rise to a local free sheaf $\sV_k$ on $Y$ with the
Gau\ss-Manin connection
$$
\nabla: \sV_k \>>> \sV_k \otimes\Omega^1_Y(\log(S+T)),
$$
where we assume that $\sV_k$ is the quasi-canonical extension of
$$
(R^kg_*\C_{Z_0})\otimes_\C \sO_{Y\setminus (S\cup T)},
$$
i.e. that the real part of the eigenvalues of the
residues around the components of $S+T$ lie in $[0,1)$.

By \cite{Sch} $\sV_k$ carries a filtration $\sF^p$ by coherent
subsheaves. If the monodromies around the components of $S+T$
are not unipotent the $\sF^p$ are not necessarily subbundles.
However this is the case outside of the singular locus of $S+T$.
By abuse of notations, we will drop the assumption that $Y$ is
projective in the first part of this section, leave out a
codimension two subscheme $W$ and assume that
$f$, $f'$ and $g$ are flat and that $S+T$ is non-singular.

So the induced graded sheaves $E^{p,k-p}$ are locally free,
and they carry a Higgs structure with logarithmic poles along
$S+T$. Let us denote it by
$$
({\mathfrak g \mathfrak r}_{\sF}(\sV_k), {\mathfrak g \mathfrak
r}_{\sF}(\nabla))=(E,\theta)=\left(
\bigoplus_{q=0}^{k}E^{k-q,q}\ , \
\bigoplus_{q=0}^{k}\theta_{k-q,q}\right).
$$
As in \cite{V-Z} we will consider a second
system of sheaves related to $Z$ and to the pair $(X,H)$.
We define
$$
F^{p,q}= R^qf'_*(\delta^*(\Omega_{X/Y}^p(\log \Delta))\otimes
{\sL'}^{(-1)})/_{\rm torsion}.
$$
Here, for $\eta=0, \ldots , \nu-1$,
the invertible sheaves ${\sL'}^{(-\eta)}$ are defined as
$$
{\sL'}^{(-\eta)}= {\sL'}^{-\eta}\otimes
\sO_{X'}\left(\left[\frac{\eta\cdot(H'+E')}{\nu}\right]\right)=
{\sL'}^{-\eta}\otimes
\sO_{X'}\left(\left[\frac{\eta\cdot E'}{\nu}\right]\right).
$$
As well-known (see for example \cite{Gri}, page 130) the bundles
$E^{p,q}$ have a similar description:
$$
E^{p,q}=R^qg_*\Omega_{Z/Y}^p(\log \Pi).
$$
Let
\begin{gather*}
\tau_{p,q}:F^{p,q} \>>> F^{p-1,q+1}\otimes \Omega^1_Y(\log T)
\mbox{ \ \ \ and \ \ \ }\\
\tilde{\theta}_{p,q}:E^{p,q} \>>> E^{p-1,q+1}\otimes \Omega^1_Y(\log (S+T))
\end{gather*}
be the edge morphisms of the tautological exact sequences
\begin{multline}\label{exact1}
0\to {f'}^*\Omega^1_Y(\log T)\otimes
\delta^*(\Omega^{p-1}_{X/Y}(\log \Delta))\otimes {\sL'}^{(-1)}
\to \\
\delta^*({\mathfrak g \mathfrak r}(\Omega_X^p(\log
\Delta)))\otimes {\sL'}^{(-1)} \to
\delta^*(\Omega_{X/Y}^p(\log \Delta))\otimes {\sL'}^{(-1)} \to 0,
\end{multline}
and
\begin{multline}\label{exact2}
0\to {g}^*\Omega^1_Y(\log (S+T))\otimes
\Omega^{p-1}_{Z/Y}(\log \Pi)\to
{\mathfrak g \mathfrak r}(\Omega_Z^p(\log \Pi)) \to
\Omega_{Z/Y}^p(\log \Pi) \to 0,
\end{multline}
respectively, where
\begin{gather*}
{\mathfrak g \mathfrak r}(\Omega_X^p(\log \Delta))=
\Omega_X^p(\log \Delta)
/f^*\Omega^2_Y(\log T)\otimes \Omega^{p-2}_{X/Y}(\log
\Delta)\mbox{, \ \ \ and}\\
{\mathfrak g \mathfrak r}(\Omega_Z^p(\log \Pi))=
\Omega_Z^p(\log \Pi)
/g^*\Omega^2_Y(\log S+T)\otimes \Omega^{p-2}_{Z/Y}(\log \Pi).
\end{gather*}

The Gau\ss-Manin connection is the edge morphism
of
\begin{multline*}
0\to {g}^*\Omega^1_Y(\log (S+T))\otimes
\Omega^{\bullet-1}_{Z/Y}(\log \Pi)\to
{\mathfrak g \mathfrak r}(\Omega_Z^{\bullet}(\log \Pi)) \to
\Omega_{Z/Y}^{\bullet}(\log \Pi) \to 0,
\end{multline*}
hence $\theta_{p,q}=\tilde{\theta}_{p,q}$.
\begin{lemma}\label{coverings}
Let $\bullet$ stand either for ${\rm Spec}(\C)$ or for $Y$. Then
the group $\Z/\nu$ acts on $\psi_*\Omega^p_{Z/\bullet}(\log
(\Pi+ \psi^* H'))$ and on $\psi_*\Omega^p_{Z/\bullet}(\log (\Pi))$.
One has a decomposition in sheaves of eigenvectors
\begin{gather*}
\psi_*\Omega^p_{Z/\bullet}(\log (\Pi+ \psi^* H')) \cong
\bigoplus_{\eta=0}^{\nu-1} \Omega^p_{X'/\bullet}(\log
(\Delta'+\Sigma'+ H'))\otimes {\sL'}^{(-\eta)}
\mbox{ \ \ \ and}\\
\psi_*\Omega^p_{Z/\bullet}(\log  \Pi ) \cong
\Omega^p_{X'/\bullet}(\log (\Delta'+\Sigma')) \oplus
\bigoplus_{\eta=1}^{\nu-1} \Omega^p_{X'/\bullet}(\log
(\Delta'+\Sigma'+ H'))\otimes {\sL'}^{(-\eta)},
\end{gather*}
compatible with the tautological sequences.
\end{lemma}
\begin{proof}
By \cite{E-V}, 3.21 and 3.22, there are natural inclusions
$$
\psi^*\Omega^p_{X'/\bullet}(\log(\Delta'+\Sigma'+H')) \>>>
\Omega_{Z/\bullet}^p(\log(\Pi+\psi^*H')),
$$
and $R^\beta\psi_*\Omega_{Z/\bullet}^p(\log(\Pi+\psi^*H'))=0$,
for $\beta>0$. In fact, in \cite{E-V} this is just stated for
$\bullet={\rm Spec}(\C)$, but the general case follows by induction,
considering the tautological sequences. Since $\Z/\nu$ acts on
$\psi_*\sO_Z$ with
$$
\psi_*\sO_Z=\bigoplus_{\eta=0}^{\nu-1} {\sL'}^{(-\eta)}
$$
as decomposition in sheaves of eigenvectors, one obtains the
first decomposition in \ref{coverings}. $H'$ is totally ramified
in $Z$. Hence there is an exact sequence
$$
0\to \psi_*\Omega^p_{Z/\bullet}(\log \Pi)\to
\psi_*\Omega^p_{Z/\bullet}(\log (\Pi+ \psi^* H')) \to
\Omega^{p-1}_{H'/\bullet}(\log (\Delta'+\Sigma')|_{H'})
$$
and the two sheaves on the right hand side differ only in the
$\Z/\nu$ invariant part.
\end{proof}
\begin{lemma}\label{6.1} Using the notations introduced above,
let
$$
\iota:\Omega^1_Y(\log T)\>>>\Omega^1_Y(\log (S+T))
$$
be the natural inclusion. Then there exist morphisms
$\rho_{p,q}: F^{p,q} \to E^{p,q}$
such that:
\begin{enumerate}
\item[i)] The diagram
$$
\begin{CD}
E^{p,q} \> \theta_{p,q} >> E^{p-1, q+1} \otimes \Omega^{1}_{Y}
(\log (S+T)) \\
\A \rho_{p,q} AA \A A \rho_{p-1,q+1} \otimes \iota A \\
F^{p,q} \> \tau_{p,q} >> F^{p-1,q+1} \otimes \Omega^{1}_{Y} (\log T).
\end{CD}
$$
commutes.
\item[ii)] There is an invertible sheaf $\sA$, semi-ample and
ample with respect to $Y\setminus T$, an effective
divisor $D'$, and an injection $\sA(D')\to F^{n,0}$, which is
an isomorphism over $Y\setminus T$.
\item[iii)] $\tau_{n,0}$ induces a morphism
$$
\tau^\vee:T_Y( - \log T)=(\Omega_Y^1(\log T))^\vee \>>>
F^{{n,0}^\vee}\otimes F^{n-1,1},
$$
which coincides over $Y\setminus (S\cup T)$ with the Kodaira-Spencer map
$$
T_Y( - \log T) \>>> R^1f_*T_{X/Y}(-\log \Delta).
$$
In particular this morphism is injective.
\item[iv)] The morphisms $\rho_{n-m,m}$ are injective, for all $m$.
\item[v)] \hfill
$\displaystyle
\left( \bigoplus_{q=0}^{n}E^{n-q,q}\ , \
\bigoplus_{q=0}^{n}\theta_{n-q,q}\right)
$\hspace*{\fill}\\
is a Higgs bundle with logarithmic poles along $S+T$, induced by a variation
of Hodge structures with finite monodromy around the components of $S$.
\end{enumerate}
\end{lemma}
\begin{remark} \label{6.1b}
Instead of \ref{6.1}, iii) and iv) we will use
later just the injectivity of $\tau^\vee$ and of
$\rho_{n-m,m}$ for $m=0$ and $m=1$.
\end{remark}
\begin{proof} The proof is similar to the proof of 3.2 in
\cite{V-Z}. It is well-known, that the bundle in v) is the Higgs
bundle for the variation of Hodge structures on
$R^ng_*\C_{Z_0}$. The condition on the monodromy
follows from \ref{5.5}. By \ref{coverings} the sheaf
$$
R^q{f'}_*(\Omega_{X'/Y}^p(\log
(H'+\Delta'+\Sigma'))\otimes {\sL'}^{(-1)})
$$
is a direct factor of $E^{p,q}$. The morphism $\rho_{p,q}$
is induced by the natural inclusions
\begin{multline}\label{inclusion}
\delta^*\Omega_{X/Y}^p(\log \Delta)\to
\delta^*\Omega_{X/Y}^p(\log (\Delta+\Sigma))\\
\to \Omega_{X'/Y}^p(\log (\Delta'+\Sigma'))\to
\Omega_{X'/Y}^p(\log (H'+\Delta'+\Sigma')).
\end{multline}
Over $Y\setminus (S\cup T)$ the kernel of $\rho_{n-m,m}$ is a quotient of
the sheaf
$$
R^{m-1}(f'|_{H'})_*(\Omega_{H'/Y}^{n-m-1} \otimes {\sL'}^{-1}|_{H'}).
$$
Since the relative dimension of $H'$ over $Y$ is $n-1$ and since $\sL'$ is
fibrewise ample, the latter is zero by the
Akizuki-Kodaira-Nakano vanishing theorem. So $\rho_{n-m,m}$ is
injective, as claimed in iv).

The injective morphism in (\ref{inclusion}) also exist for $Y$ replaced
by ${\rm Spec}(\C)$, and the exact sequence (\ref{exact1}) is a
subsequence of
\begin{multline}\label{exact4}
0\to {f'}^*\Omega^1_Y(\log (S+T))\otimes
\Omega^{p-1}_{X'/Y}(\log (H'+ \Delta +\Sigma))\otimes
{\sL'}^{(-1)} \to \\
{\mathfrak g \mathfrak r}(\Omega_{X'}^p(\log (H'+ \Delta+\Sigma
)))\otimes {\sL'}^{(-1)} \to
\Omega_{X'/Y}^p(\log(H'+ \Delta+\Sigma))\otimes {\sL'}^{(-1)} \to 0.
\end{multline}
Finally by \ref{coverings} this sequence is obtained by taking
the sheaves of eigenvectors in the direct image of the
exact sequence (\ref{exact2}) under $\psi:Z\to X'$. One obtains
i).

By definition $F^{n,0}= f'_*(\delta^*(\Omega^n_{X/Y}(\log \Delta))\otimes
{\sL'}^{(-1)})$. Comparing the first Chern classes for the
tautological sequence for $f$ one finds
$$
F^{n,0}= f'_*(\delta^*(\omega^n_{X/Y}(\Delta_{\rm red} - \Delta))\otimes
{\sL'}^{(-1)}).
$$
Recall that $f$ is smooth over $Y\setminus B$, for the divisor $B$
considered in \ref{4.3}, c). Hence
$$f^*B \geq - \Delta_{\rm red} + \Delta$$
and $\Omega^n_{X/Y}(\log \Delta)$ contains
$\omega_{X/Y}(-f^*B)$.
Moreover, by \ref{4.3}, c), $D'=D-B$ is effective.
By definition,
$\sL = \omega_{X/Y} \otimes f^* \sA(D'+B)^{-1}$, and
$$
{\sL'}^{(-1)} = {\sL'}^{-1}\otimes
\sO_{X'}\left(\left[\frac{E'}{\nu}\right]\right).
$$
Therefore $\delta^*(\Omega^n_{X/Y}(\log \Delta))\otimes
{\sL'}^{(-1)}$ contains $\omega_{X/Y}(-f^*B)\otimes
{\sL'}^{(-1)}$, hence the sheaf
$$
{f'}^*(\sA(D'))\otimes
\sO_{X'}\left(\left[\frac{E'}{\nu}\right]\right),
$$
and ii) holds true. For iii), recall that over $Y\setminus (S\cup T)$
the sheaf ${\sL'}^{(-1)}$ is nothing but
$$
{\sL'}^{-1}=\delta^*(\sL^{-1}).
$$
Since $R^\mu\delta_*\sO_{X'}=0$, by the projection formula
the morphism
$$
(\tau_{n,0}\otimes {\rm id}_{\sA(D')^{-1}})|_{Y\setminus (S\cup T)}
$$
is the restriction of the edge morphism of the short exact sequence
\begin{equation*}
0\to {f}^*\Omega^1_U \otimes \Omega^{n-1}_{V/U}\otimes \omega_{V/U}^{-1}\to
{\mathfrak g \mathfrak   r}(\Omega_{V}^n)\otimes {\omega_{V/U}}^{-1} \to
\Omega_{V/U}^{n} \otimes \omega_{V/U}^{-1} \to 0.
\end{equation*}
Since $f|_V$ is smooth with $n$-dimensional fibres, the sheaf on
the right hand side is $\sO_V$ and the one on the left hand side
is $f^*\Omega^1_U\otimes T_{V/U}$.
Tensoring with
$$
f^*T_U=f^*(\Omega_U^{r-1} \otimes \omega_U^{-1})
$$
and dividing by the kernel of the wedge product
$$
f^*\Omega_U^1\otimes f^*(\Omega_U^{r-1} \otimes \omega_U^{-1})
\>>> \sO_V
$$
on the left hand side, one obtains an exact
sequence
\begin{equation}\label{taut1}
0 \>>> T_{V/U} \>>> \sG \>>> f^*T_U  \>>> 0,
\end{equation}
where $\sG$ is a quotient of
${\mathfrak g \mathfrak r}(\Omega_{V}^n)\otimes
{\omega_{V}}^{-1} \otimes f^* \Omega^{r-1}_U$.
By definition, the restriction to $Y\setminus (S\cup T)$ of
the morphism considered in iii) is the first edge
morphism in the long exact sequence, obtained by applying
$R^\bullet f_*$ to (\ref{taut1}).

The wedge product induces a morphism
$$
\Omega_{V}^n \otimes {\omega_{V}}^{-1} \otimes f^*
\Omega^{r-1}_U \>>> \Omega_V^{n+r-1} \otimes \omega_V^{-1}=T_V.
$$
Since $r=\dim(U)$ this morphism factors through $\sG$.
Hence the exact sequence (\ref{taut1}) is isomorphic
to the tautological sequence
\begin{equation}\label{taut2}
0\>>> T_{V/U} \>>> T_V \>>> f^*T_U \>>> 0.
\end{equation}
The edge morphism $T_U \to R^1f_*T_{V/U}$ of (\ref{taut2})
is the Kodaira-Spencer map.
Since we assumed $U$ to be generically finite over the moduli
space, this morphism is injective.
\end{proof}
Let us return to the case ``$Y$ projective''. We will choose for
$E^{p,q}$ and $F^{p,q}$ the maximal coherent extension
of the sheaves defined above outside of a codimension two
subvariety of $Y$. Of course, the morphisms $\theta_{p,q}$,
$\tau_{p,q}$ and $\rho_{p,q}$ extend, and the properties i) -
v) in \ref{6.1} remain true.

By \cite{Sim}, page 12, $\theta\wedge\theta=0$ hence the image
of the composite
$$
\theta_{n-q+1,q-1}\circ \cdots \circ \theta_{n,0}:
E^{n,0} \>>> E^{n-q,q}\otimes \bigotimes^{q} \Omega^1_Y(\log
(S+T))
$$
factors through
$$
\theta^q: E^{n,0} \>>> E^{n-q,q}\otimes S^{q} \Omega^1_Y(\log
(S+T)).
$$
By \ref{6.1}, ii) $\sA(D')$ is a subsheaf of $F^{n,0}$ and hence
of $E^{n,0}$, and one obtains a morphism
\begin{multline*}
\sA(D') \>>> \rho_{n-q,q}(F^{n-q,q})\otimes S^{q} \Omega^1_Y(\log T)
\> \subset >>  E^{n-q,q}\otimes S^{q} \Omega^1_Y(\log T)\\
\> S^m(\iota) >> E^{n-q,q}\otimes S^{q} \Omega^1_Y(\log (S+T)),
\end{multline*}
and thereby a morphism
$$
{\tau'}^q: S^{q}(T_Y(-\log T)) \>>> E^{n-q,q}\otimes \sA(D')^{-1}.
$$
The pullback of ${\tau'}^q$, via $\gamma:\C \>>> Y\setminus T \>>> Y$,
composed with the $q$-th tensor power of the differential of $\gamma$
$$
d\gamma^{q}:T_\C^{q} \>>> \gamma^*(S^{q} T_Y(-\log T))
$$
gives
$$
\tilde{\tau}^q: T_\C^{q} \>>> \gamma^*(
E^{n-q,q}\otimes \sA(D')^{-1}).
$$
We choose
$$
m={\rm Min}\{q\in \N; \ \tilde{\tau}^{q+1}(T_\C^{
q+1})=0\}
$$
and put $\tau={\tau'}^m$ and $\tilde{\tau}=\tilde{\tau}^m$.

The morphism ${\tau'}^1$ factors like
$$
T_Y(-\log T) \>>> F^{n-1,1}\otimes \sA(D')^{-1}
\>  \rho_{n-1,1} >> E^{n-1,1}\otimes \sA(D')^{-1}.
$$
By \ref{6.1}, iii) the first of those morphisms is injective,
and by \ref{6.1}, iv) the second one as well. Therefore
${\tau'}^1$ is injective. Since we assumed $\gamma(\C)$ to be
dense, the pullback of an injective morphism of sheaves under
$\gamma$ remains injective. Hence $\tilde{\tau}^1$ is injective,
and $m>0$.

Altogether, starting from the morphism in \ref{4.4} and from a
holomorphic map $\gamma:\C \to U$ with dense image, we
constructed divisors $S$ and $T$ with the properties stated in
\ref{5.4}, and we constructed Higgs bundles which satisfy the
properties a) - d) given below.
\begin{lemma}\label{6.2}
For some $m>0$ there exist an invertible sheaf $\sA$, an
effective divisor $D'$ and a morphism of sheaves
$$
{\tau}: S^{m}T_Y(-\log T) \>>> E^{n-m,m}\otimes \sA(D')^{-1}\>>>
E^{n-m,m}\otimes \sA^{-1},
$$
such that the composite
$$
\tilde{\tau}=\gamma^*{\tau}\circ d\gamma^{m}:
T_\C^{m} \>>> \gamma^*(S^{m}T_Y(-\log T))
\>>>\gamma^*(E^{n-m,m}\otimes \sA^{-1})
$$
satisfies:
\begin{enumerate}
\item[a)] $\tilde{\tau}$ is injective.
\item[b)] \hfill $\tilde{\tau}(T_\C^{m})\subset
\sN\otimes\gamma^*(\sA^{-1})$ \hspace*{\fill}\\[.1cm]
for a sub-linebundle $\sN$ of
$$
{\rm Ker}(\gamma^*(\theta_{n-m,m}):\gamma^*(E^{n-m,m})\>>>
\gamma^*(E^{n-m-1,m+1})\otimes\Omega_\C^1(\log \gamma^{-1}(S))).
$$
\item[c)]
$$
(E,\theta)=\left( \bigoplus_{p+q=n}E^{p,q}, \theta_{p,q}\right)
$$
is the Higgs bundle, corresponding to the quasi-canonical extension
$\sV$ of $\V_0\otimes_\C\sO_{Y\setminus(S\cup T)}$ for a
geometric variation of Hodge structures $\V_0$, with finite monodromies
around the components of $S$.
\item[d)] $\gamma(\C)$ does not meet $T$.
\item[e)] $\sA$ is ample.
\end{enumerate}
\end{lemma}
At the end of the next section we will show, that those
properties lead to a contradiction to the Ahlfors-Schwarz lemma.
Hence the holomorphic map $\gamma$ can not exist.
\begin{proof}
All properties hold true, for the Higgs bundles constructed
above, with $\sA$ semi-ample and big. Choose some $\eta
>0$ such that $\sA^\eta$ contains an ample invertible sheaf
$\sA'$, and consider the Higgs bundles
$$
(E',\theta')=(E^{\otimes \eta},\theta') \mbox{ \ \ \  and \ \ \ }
(F',\tau') = (F^{\otimes \eta},\tau').
$$
Again we first consider them on $Y-W$, where $W$ is the
singular locus of $S\cup T$ and then we take the maximal
extension to $Y$. By \cite{Sim}, page 70, the morphism
$$
\theta':E^{\otimes \eta} \>>> E^{\otimes \eta}\otimes
\Omega_Y(\log (S+T))
$$
is given by
$$
\theta'= \theta\otimes {\rm id}_{E} \otimes \cdots \otimes {\rm
id}_{E} + {\rm id}_{E}\otimes \theta \otimes \cdots \otimes {\rm
id}_{E} + \cdots + {\rm id}_{E} \otimes \cdots \otimes {\rm
id}_{E}\otimes \theta,
$$
and similarly for $F'$ and $\tau'$.
The decomposition as a direct sum is
$$
\bigoplus_{p+q=k}{E'}^{p,q} \mbox{ \ \ \ and \ \ \ }
\bigoplus_{p+q=k}{F'}^{p,q},
$$
with
$$
{E'}^{p,q}=\bigoplus
\bigotimes_{i=1}^\eta E^{p_i,q_i}\mbox{ \ \ \ and \ \ \ }
{F'}^{p,q}=\bigoplus \bigotimes_{i=1}^\eta F^{p_i,q_i},
$$
where the direct sums are taken over all $p_1, \ldots , p_\eta,
q_1, \ldots , q_\eta$ with
$$
\sum_{i=1}^\eta p_i=p, \mbox{ \ \ \ and \ \ \ } \sum_{i=1}^\eta q_i = q.
$$
Again we have morphisms
$$
\rho'_{p,q}=\bigoplus\bigotimes_{i=1}^{\eta} \rho_{p_i,q_i}:
{F'}^{p,q} \>>> {E'}^{p,q},
$$
compatible with $\theta'_{p,q}$ and $\tau'_{p,q}$.
In particular $\rho'_{n\eta,0}$ is the $\eta$-th tensor product of
$\rho_{n,0}$, hence injective. The same holds true for
$\rho'_{n\eta-1,1}$ which is the direct
sum of morphisms of the form
$$
\rho_{n,0}\otimes \cdots \otimes \rho_{n-1,1} \otimes \cdots
\otimes \rho_{n,0}.
$$
The properties i) and v) in \ref{6.1} remain true, with $E$ and
$F$ replaced by $E'$ and $F'$, for $n\eta$ instead of $n$ in v). In
ii) one has an injection
$$
\sA' \>>> \sA(D')^\eta \>>> {F'}^{n\eta,0}=(F^{n,0})^{\otimes \eta}.
$$
The morphism
\begin{multline*}
\tau'_{n\eta,0}: {F^{n,0}}^{\otimes \eta} \>>>
(F^{n-1,1}\otimes F^{n,0} \otimes \cdots \otimes F^{n,0}
\oplus F^{n,0}\otimes F^{n-1,1} \otimes \cdots \otimes F^{n,0}
\oplus \\
\cdots \oplus F^{n,0}\otimes F^{n,0} \otimes \cdots
\otimes F^{n-1,1})\otimes \Omega^1_Y(\log T)
\end{multline*}
is a direct sum of morphisms of the form
$$
{\rm id}_{F^{n,0}} \otimes \cdots \otimes \theta_{n,0} \otimes
\cdots \otimes {\rm id}_{F^{n,0}},
$$
hence it induces the diagonal morphism
$$
\oplus {\tau}^\vee:T_Y(-\log T) \>>> {{F'}^{n\eta,0}}^\vee \otimes
{F'}^{n-1,1} = \bigoplus^\eta {F^{n,0}}^\vee \otimes F^{n-1,1}.
$$
In particular the injectivity of the morphisms in \ref{6.1},
iii), carries over.

As remarked in \ref{6.1b} the injectivity of $\oplus
{\tau}^\vee$,  $\rho'_{n\eta,0}$ and of $\rho'_{n\eta-1,1}$
is sufficient to perform the constructions with $E'$ and $F'$
instead of $E$ and $F$, and to obtain some $m>0$ and the
morphisms $\tau$ and $\tilde{\tau}$ satisfying the
properties a) and b), with $\sA$ replaced by $\sA^\eta$. The latter
contains the ample sheaf $\sA'$, hence e) holds true.

Finally the Higgs bundle $(E',\theta')$ comes from the locally
free extension $\sV'=\sV^{\otimes \eta}$ of $\V_0^{\otimes
\eta}\otimes_\C \sO_{Y\setminus (S\cup T)}$. The eigenvalues of the
residues of the induced connection lie in $\R_{\geq 0}$, hence
$\sV'$ is contained in the quasi-canonical extension $\sV''$.
Replacing $\sV'$ by $\sV''$, we enlarge the sheaves ${E'}^{p,q}$,
which is allowed without changing the properties a) and b).
\end{proof}
%%%%%%%%%%%%%%%%%%%%%%%% Sect. 7 %%%%%%%%%%%%%%%%%%%%
\section{Curvature estimates and the Ahlfors-Schwarz lemma}

Let $T$ be  the normal crossing divisor in \ref{6.2},
and let $T=\sum_{i=1}^\ell T_i$ be its decomposition in irreducible
components. Let $s_i$ be the section of $\sL_i=\sO_Y(T_i)$ with
zero set $T_i$.  We choose  a  hermitian metric $g_i$ on  $\sL_i$
and define
$$
r_i=-\log||s_i||_{g_i}^2 \mbox{ \ \ and \ \ }
r=r_1\cdot\cdots\cdot r_\ell.
$$
Given any constant $c>1$, by rescaling
the sections $s_i$, i.e. by replacing $s_i$ by $\epsilon\cdot
s_i$, for $\epsilon$ sufficiently small, one may assume that
$r_i \geq c$.

On the ample invertible linebundle $\sA$ in \ref{6.2} we
choose a metric $g$ such that the curvature form $\Theta(\sA,g)$
is positive definite. For a positive integer $\alpha$ we define a
new metric $ g_{\alpha}=g\cdot r^{\alpha} $ on $\sA|_{Y\setminus
T}.$

Recall that a hermitian form $\omega_\alpha$ on $T_Y(-\log T)$ is
continuous and positive definite, if each point in $Y$ has a
neighborhood $U$ with local coordinates $z_1,\ldots ,z_n$, such
that $T\cap U$ is the zero set of $z_1\cdot\cdots\cdot z_k$ and
such that, writing
\begin{gather*}
\iota_1=\cdots =\iota_k=1
\mbox{ \ \ \ and \ \ \ } \iota_{k+1}=\cdots \iota_n=0,
\\
\omega_{\alpha}=\sqrt{-1}\sum_{1\leq  i,j\leq n }a_{i,j}\frac
{dz_i}{z_i^{\iota_i}}\wedge \frac{d\bar z_j}{\bar{z_j}^{\iota_j}}
\end{gather*}
for a continuous and positive definite hermitian matrix
$(a_{ij})_{1\leq i,j\leq n}$.
\begin{lemma}\label{7.1} Rescaling the $s_i$, if necessary,
there exists a continuous and positive
definite hermitian  form  $\omega_{\alpha}$ on $T_Y(-\log T)$
with
$$
r^2\Theta(\sA|_{Y\setminus T}, g_{\alpha})\geq \omega_{\alpha}.
$$
\end{lemma}
\begin{proof} We recall the formula for the curvature
calculation of a linebundle with a metric $(\sL,g)$ (see
for example \cite{Dem}, 7.1). Let
$$
\sL|_U\simeq U\times \C
$$
be a local trivialization of $\sL$ and let $s_U $ be a holomorphic
section of $\sL|_U,$ which does not vanish in any point of $U$.
Then $s_U$ corresponds to a holomorphic function $h_U$ on $U,$
and the metric $g$ is given by
$$
||s_u||^2_g=|h_U|^2e^{-\phi}.
$$
The curvature  $\Theta(\sL,g)$ is given by
$$
\Theta(\sL,g)=\frac {\sqrt{-1}}{2\pi}\partial\bar\partial \phi.
$$
Applying this formula (see also \cite{Lu}, proof of 3.1), one finds
\begin{gather*}
\Theta(\sA,g_{\alpha})=\Theta(\sA,gr^{\alpha})=\Theta(\sA,ge^{-(-\alpha\log
r)})=\Theta(A,g)-\frac{\sqrt{-1}\alpha}{2\pi} \partial\bar\partial\log r\\
=\Theta(\sA,g)
-\sum_{i=1}^\ell \frac{\sqrt{-1}\alpha}{2\pi}
\partial\bar\partial\log r_i
=\Theta(\sA,g)-\sum_{i=1}^\ell\frac{\sqrt{-1}\alpha}{2\pi}\partial\frac{\bar\partial
r_i}{r_i}\\
=\Theta(\sA,g)-\sum_{i=1}^\ell
\frac{\alpha\Theta(L_i,g_i)}{r_i}+ \frac{\sqrt{-1}\alpha}{2\pi}
\frac{\partial r_i\wedge \bar\partial r_i}{r_i^2}.
\end{gather*}
Rescaling the sections $s_i$ one can assume that the $r_i$ are
larger than a large constant $c>1$, hence that
$$
\Theta(\sA,g)-\sum_{i=1}^\ell
\omega'_{\alpha} := \frac{\alpha\Theta(L_i,g_i)}{r_i}
$$
is a continuous and positive definite (1,1)-form on $Y.$ Moreover
$$
\Theta(\sA,g_{\alpha})=\omega'_{\alpha}+\sum_{i=1}^\ell
\frac{\sqrt{-1}\alpha}{2\pi}
\frac{\partial r_i\wedge \bar\partial r_i}{r_i^2}\geq
\omega'_\alpha+\frac{\sqrt{-1}\alpha}{2\pi}\sum_{i=1}^\ell
\frac{\partial r_i\wedge \bar\partial r_i}{r^2}.
$$
The (1,1) form
$$
\frac{\sqrt{-1}\alpha}{2\pi}\sum_{i=1}^\ell  \partial r_i\wedge
\bar\partial r_i
$$
is clearly positive semi-definite on $Y\setminus T$.

\begin{claim} \label{7.2} Assume again that $T\cap U$ is the
zero set of $z_1\cdot \cdots \cdot z_k$ for local coordinates
$z_1,\ldots , z_n$ on $U$. Then in a
small neighborhood of $T\cap U$ the form
$$
\frac{\sqrt{-1}\alpha}{2\pi}\sum_{i=1}^\ell  \partial r_i\wedge
\bar\partial r_i
$$
is positive definite on the subspace of $T_Y(-\log T)|_U$
spanned by
$$
\{z_1 \partial_{z_1},\dots ,z_k\partial_{z_k}\}.
$$
\end{claim}
\begin{proof} Near  $T_i$ the section $s_i$ can be expressed as
$$
s_i=z_it_i,\quad ||s_i||^2_{g_i}=z_i\bar z_i ||
t_i||^2_{g_i}=z_i\bar z_i f_i,
$$
where $t_i$ a local basis of  $L_i$  and where $f_i$ is a positive
Function. So,
\begin{gather*}
r_i=-\log||s_i||^2_{g_i}=-\log z_i-\log \bar z_i -\log f_i,\\
\partial
r_i=-\frac{dz_i}{z_i}-\frac{1}{f_i}\sum_{j=1}^n\frac{\partial
f_i}{\partial z_j}dz_j,
\end{gather*}
and
$$
\bar \partial r_i=-\frac{d\bar z_i}{\bar
z_i}-\frac{1}{f_i}\sum_{j=1}^n\frac{\partial f_i}{\partial \bar
z_j}d\bar z_j.
$$
So the leading term in
$$
\frac{\sqrt{-1}\alpha}{2\pi}\sum_{i=1}^\ell  \partial r_i\wedge
\bar\partial r_i
$$
near $T\cap U$ is
$$
\frac{\sqrt{-1}\alpha}{2\pi}\sum_{i=1}^k \frac{dz_i}{z_i}\wedge
\frac{d\bar z_i}{\bar z_i}.
$$
Obviously this form is positive definite on the subspace spanned
by
$$
\{z_1\partial_{z_1},\ldots ,z_k\partial_{z_k}\}.
$$
\end{proof}
Since we assumed that $r\geq 1$,
$$
r^2\Theta(A|_{Y\setminus T}, g_{\alpha})\geq
r^2\omega'_{\alpha}+\frac{\sqrt{-1}\alpha}{2\pi}\sum_{i=1}^\ell  \partial
r_i\wedge\bar\partial r_i\geq
\omega'_{\alpha}+\frac{\sqrt{-1}\alpha}{2\pi}\sum_{i=1}^\ell
\partial r_i\wedge\bar\partial r_i.$$
By \ref{7.2} the (1,1) form
$$
\omega_\alpha=\omega'_{\alpha}+\frac{\sqrt{-1}\alpha}{2\pi}\sum_{i=1}^\ell
\partial r_i\wedge\bar\partial r_i
$$
is continuous and positive definite on $T_Y(-\log T).$
\end{proof}
Let $\gamma: \C\to Y\setminus T$ be the holomorphic map
with Zariski dense image, considered in \ref{6.2} and let $t$ be
the global coordinate on $\C.$ We take the ample bundle $\sA$ on
$Y$ with the metric $g_{\alpha}$ on $Y\setminus T$ and the
hermitian metric $\omega_{\alpha}$ on $T_Y(-\log T)$ from
lemma \ref{7.1}. Writing again
$$
d \gamma: T_\C \to \gamma^* T_{Y}(-\log T)
$$
for the differential, one finds
$$
\gamma^*\omega_{\alpha}=\sqrt{-1}||d \gamma
(\partial_t)||^2_{\gamma^*\omega_{\alpha}}dt\wedge d{\bar t},
$$
and \ref{7.1} implies:
\begin{corollary}\label{7.3} \ \
$\gamma^*r^2\Theta(\sA|_{Y\setminus T},g_{\alpha})\geq
\sqrt{-1}||d \gamma
(\partial_t)||^2_{\gamma^*\omega_{\alpha}}dt\wedge d{\bar t}.$
\end{corollary}
Let us return to the morphism of sheaves in \ref{6.2}:
\begin{gather*}
\tau: S^mT_Y(-\log T)\>>> E^{n-m, m}\otimes
(\sA(D'))^{-1}\hookrightarrow E^{n-m, m}\otimes \sA^{-1} \\
\tilde \tau:=\gamma^*\tau\circ(d\gamma)^{m}: T^{
m}_{\C}\>>> \gamma^*S^m T_{Y}(-\log T)\hookrightarrow
\gamma^*(E^{n-m,m}\otimes \sA^{-1}).
\end{gather*}
By \ref{6.2}, c)  $E^{n-m,m}$ is a sub quotient of
the quasi-canonical extension of a geometric variation of Hodge
structures $\V_0$ on $Y\setminus S\cup T$. By Kawamata's
construction (see \cite{Vie}, 2.5) one finds a finite morphism
$ \pi: Y'\to Y$ with $Y'$ non-singular and $S'+T'=\pi^*(S+T)$ a
normal crossing divisors such that the local monodromies of the
pullback $\pi^*\V_0$ around $S'+T'$ are unipotent.
For the discriminant $\Delta(Y'/Y)$ of $\pi:Y'\to Y$ both,
$$
\Delta(Y'/Y)+S+T \mbox{ \ \ \ and \ \ \ }\pi^*(\Delta(Y'/Y)+S+T)
$$
are normal crossing divisors. Moreover, for a component $T_i$ of
$S+T$ there exists some $\mu_i$ with
$$
\pi^*T_i=\mu_i\cdot (\pi^*T_i)_{\rm red}.
$$
Since we assumed the local monodromy of $\V_0$ around the
components of $S$ to be of finite order, the local monodromy
of $\pi^*\V_0$ around the components of $S'=\pi^*(S)$ is
trivial, hence $\pi^*\V_0$ extends to a variation
of Hodge structures $\V'_0$ across $S'.$ Let $h$ and $h'$ denote
the Hodge metrics on $\V_0$ and $\V'_0$, respectively. We use
the same notation for the induced metric on the Higgs bundles
$\bigoplus {E}^{p,q}$ and $\bigoplus {E'}^{p,q}$,
where the latter is given by sub quotients of the
canonical extension of $\V'_0$ across $T'=\pi^*T.$
We have an inclusion of sheaves
$$
{\iota}: (\pi^*E^{n-m,m}, \pi^*h) \hookrightarrow ({E'}^{n-m,m},
h'),
$$
such that $\pi^*(h)={\iota}^*(h')$ on $Y'\setminus S'\cup T'$.

Consider the diagram of morphisms of analytic spaces
\begin{equation}\label{diag1}
\begin{CD}
\C'\> {\gamma'}>> Y'\\
\V \pi' VV \V \pi VV \\
\C \>{\gamma}>> Y
\end{CD}
\end{equation}
where $\C'$ is obtained as a normalization of the fibre product.
Hence if $U\subset \C$ is a sufficiently small neighborhood of
$t_0\in\gamma^{-1}(S)$, then for each $t'_0\in {\pi'}^{-1}(t_0)$
there exists a connected component $U'\subset {\pi'}^{-1}(U)$ and
a coordinate function $t'$ on $U'$, for which the map $\pi':
U'\to U$ is given by
\begin{equation}\label{local}
t-t_0=\pi'(t)=(t'-t'_0)^{\mu_0}, \mbox{ \ \ \ for some \ } \mu_0 \in \N-\{0\}.
\end{equation}
By \ref{6.2}, b), $\tilde{\tau}(T^{m}_\C)$ is contained
in an invertible linebundle $\sN\otimes \gamma^*(\sA^{-1})$,
where $\sN$ is a sub-linebundel of the kernel of
$\gamma^*(\theta_{n-m,m})$. If
$$
\theta'_{n-m,m}: {E'}^{n-m,m}\>>> {E'}^{n-m-1,m+1}\otimes
\Omega^1_{Y'}(\log \pi^*(T))
$$
denotes the Higgs structure on $Y'$, we have a commutative diagram
$$
\begin{CD}
{\gamma'}^* {E'}^{n-m,m} \> {\gamma'}^*(\theta'_{n-m,m}) >>
{\gamma'}^* {E'}^{n-m-1,m+1} \otimes \Omega^1_{\C'}(\log S')\\
\A {\gamma'}^*(\iota) A \subset A \A  {\gamma'}^*(\iota)\otimes {\pi'}^*
A \subset A \\
{\gamma'}^* \pi^*{E}^{n-m,m} \> {\gamma'}^*\pi^*(\theta_{n-m,m}) >>
{\gamma'}^* \pi^* {E}^{n-m-1,m+1} \otimes \pi^* \Omega^1_{\C}(\log S).
\end{CD}
$$
So $\iota$ induces an inclusion
$$
{\pi'}^*{\rm ker}(\gamma^*(\theta_{n-m,m})) \> \subset >>
{\rm ker }({\gamma'}^*(\theta'_{n-m,m})),
$$
hence there exists a sub-linebundle
$$
\sN'\subset {\rm ker }({\gamma'}^*(\theta'_{n-m,m}))
$$
with
\begin{equation}\label{addlocal}
{\pi'}^* \tilde{\tau}(T^{m}_\C) \subset {\pi'}^*(\sN) \otimes
{\gamma'}^*(\pi^*\sA^{-1}) \subset \sN' \otimes
{\gamma'}^*(\pi^*\sA^{-1}).
\end{equation}
As in \cite{V-Z}, 1.1, using $\theta'_{n-m,m}(\sN')=0$ and
P. Griffith's estimates for the curvature of the Hodge metric
(\cite{Gri}, chapter II) one obtains:
\begin{lemma}\label{7.11} The curvature
$\Theta(\sN',h'|_{\sN'})$ of the restricted Hodge metric on
$\sN'$ is negative semi-definite on $Y'\setminus T'$.
\end{lemma}
The Hodge metric $h$ defines a metric
$h\otimes g^{-1}_{\alpha}$ on $E^{n-m,m}\otimes
\sA^{-1}|_{Y\setminus S\cup T}.$
By \ref{6.2}, a), $\tilde\tau\not=0$ and since $\gamma(\C)$ is
Zariski dense in $Y$ we may define a non-zero positive
semi-definite K\"ahler form $\sqrt{-1}c(t)dt\wedge d\bar t $ on
$\C\setminus \gamma^{-1}(S)$ by choosing
$$
c(t)=||\tilde \tau((\partial_t)^{m})||^{\frac
{2}{m}}_{\gamma^*(h\otimes g^{-1}_{\alpha})}.
$$
\begin{lemma}\label{7.12}
Let $\mu$ denote the lowest common multiple of all the
ramification orders of components of $\pi^*(S)$ over $Y$.
Then there exists an effective divisor $\Pi$ on $\C$ (i.e. a
locally finite sum $\sum\beta_i P_i$ with $\beta_i \geq 0$) and
a linebundle $\sN^{(\mu)}$ on $\C$ with
\begin{equation*}
\tilde{\tau}(T_\C^{m})^{\mu} \otimes \sO_\C(\Pi)=
\sN^{(\mu)}\otimes \gamma^*\sA^{-\mu} \mbox{ \ \ and \ \ }
{\pi'}^*\sN^{(\mu)}={\sN'}^\mu.
\end{equation*}
\end{lemma}
\begin{proof} By (\ref{addlocal}) ${\pi'}^*\tilde{\tau}(T_{\C}^{m})$
is a subsheaf of $\sN'\otimes{\gamma'}^*(\pi^*\sA^{-1})$.
Using the description of $\pi'$ in (\ref{local}),
we choose for a given point $t'_0 \in {\pi'}^{-1}\gamma^{-1}(S)$
a small neighborhood $U'$ and some $\rho\in \N$ with
$$
{\pi'}^*\tilde{\tau}(T_{\C}^{m})|_{U'} \otimes
\sO_{U'}(\rho\cdot t'_0) = \sN'\otimes
{\gamma'}^*(\pi^*\sA^{-1})|_{U'}.
$$
The number $\frac{\rho}{\mu_0}$ is determined by the monodromy of $\V_0$
around the component of $S$ containing $\gamma(t_0)$, hence it is
independent of the point $t'_0 \in {\pi'}^{-1}(t_0)$.
Since the ramification order $\mu_0$ in (\ref{local}) divides $\mu$
we may choose $\Pi$ to be the effective divisor with $\Pi|_U=
\frac{\rho\cdot\mu}{\mu_0}\cdot t_0$ and
$$
\sN^{(\mu)}=\tilde{\tau}(T_\C^{m})^{\mu}
\otimes\sO_{\C}(\Pi)\otimes \gamma^*\sA^\mu.
$$
\end{proof}
Outside of ${\pi'}^*\Pi$ the metrics ${\gamma'}^* {h'}^\mu$
and ${\pi'}^*\gamma^* h^{\mu}$ on ${\sN'}^\mu$ coincide,
hence $\gamma^* h^{\mu}$ extends to a metric $h^{(\mu)}$ on
$\sN^{(\mu)}$ and
$$
c(t)=||\tilde \tau((\partial_t)^{m})^{\mu}||^{\frac
{2}{m\cdot \mu}}_{h^{(\mu)} \otimes \gamma^* g^{-\mu}_{\alpha}}.
$$
In particular $\sqrt{-1}c(t)dt\wedge d\bar t$ defines a
semi-definite K\"ahler form on $\C$. The induced metric $F$ is a
singular metric in the sense described in \cite{Dem},
7.1, or \cite{L-Y}, Section 2. The curvature current of
$T_{\mathbb C}$ is then defined to be  the closed (1,1)-current
$$\Theta(T_{\mathbb C}, F)=
-\frac{\sqrt{-1}}{2\pi}\partial\bar\partial\log c(t).
$$
\begin{lemma}\label{7.10}
There exists some $\epsilon' >0$ with
$$
-\Theta(T_{\mathbb C}, F)\geq
\epsilon'\gamma^*\Theta(\sA|_{Y\setminus T}, g_{\alpha})
$$
in the sense of currents.
\end{lemma}
\begin{proof}
Let $[\Pi]$ denote the current of integration over the
divisor $\Pi$. As in \cite{Dem}, proof of 7.2,
one defines a singular metric $|s|^2$ on sections of
$\sO_{\C}(\Pi)$ by taking the square of the modulus of $s$
viewed as a complex valued function. By the
Lelong-Poincar\'{e} equation $[\Pi]$ is the
curvature current of this metric. One finds
$$
\Theta(T_\C^{m\cdot \mu},F^{m\cdot\mu}) +
{\gamma}^*\Theta(\sA^\mu|_{Y\setminus T},g_{\alpha\cdot\mu}) + [\Pi] =
\Theta(\sN^{(\mu)},h^{(\mu)}).
$$
By \cite{L-Y}, section 2, the curvature current of a singular
metric on a holomorphic linebundle on a complex manifold is
compatible with pullback under holomorphic maps. Hence
$$
{\pi'}^*\Theta(\sN^{(\mu)},h^{(\mu)})=
\Theta({\sN'}^{\mu},{h'}^{\mu})=
\mu \cdot \Theta({\sN'},{h'}).
$$
By \ref{7.11} the latter is negative semi-definite, hence
$\Theta(\sN^{(\mu)},h^{(\mu)})\leq 0$.
Moreover, $[\Pi]\geq 0$ in the sense of currents, hence
$$
- \Theta(T_{\mathbb C}, F)=
- \frac{1}{m\cdot\mu}\Theta(T^{m\cdot \mu}_{\mathbb
C}, F^{m\cdot \mu})
\geq \frac{1}{m}{\gamma}^*\Theta(\sA|_{Y\setminus T},g_\alpha).
$$
\end{proof}
\begin{lemma}\label{7.9} For $\alpha\gg 1$
there exists some $\epsilon >0 $ with
$$
\gamma^*\Theta(\sA|_{Y\setminus T},g_{\alpha})\geq \epsilon
\sqrt{-1}c(t)dt\wedge d\bar t.
$$
\end{lemma}
\begin{proof} We will use the notations from \ref{7.12}, in
particular the metric $h^{(\mu)}$ on $\sN^{(\mu)}$. Recall that
$$
c(t)=||\tilde \tau((\partial_t)^{m})^{\mu}||^{\frac
{2}{m\cdot \mu}}_{h^{(\mu)} \otimes \gamma^* g^{-\mu}_{\alpha}},
$$
and that by \ref{7.3}, for all $\alpha >0$
$$
\gamma^*\Theta(\sA|_{Y\setminus T},g_\alpha) \geq
\sqrt{-1}\gamma^*r^{-2}||d\gamma(\partial_t)||^2_{\gamma^*\omega_\alpha}
dt\wedge d\bar{t}.
$$
Hence in order to show \ref{7.9} it remains to verify that
for $\alpha\gg 1,$ there exists some $\epsilon > 0$ with
\begin{multline}\label{equation1}
\gamma^*r^{-2} ||d\gamma(\partial_t)||^2_{\gamma^*\omega_{\alpha}}\geq
\epsilon \gamma^*r^{-\frac{\alpha}{m}}||\tilde
\tau((\partial_t)^{m})^{\mu}||^{\frac
{2}{m\cdot \mu}}_{h^{(\mu)}\otimes \gamma^* g^{-\mu}}\\
=\epsilon ||\tilde \tau((\partial_t)^{m})^{\mu}||^{\frac
{2}{m\cdot \mu}}_{h^{(\mu)}\otimes \gamma^* g^{-\mu}_\alpha}.
\end{multline}
Given a point $p\in Y$ choose a small polydisk $U$ with
coordinates $z_1, \ldots, z_n$, in such a way that the divisors
$T\cap U$ and $S\cap U$ are defined by the equation
$$
z_1\cdot \cdots \cdot z_{k}=0 \mbox{ \ \ and \ \ }
z_{k+1}\cdot\cdots\cdot z_{k+k'}=0.
$$
Let $\pi: Y'\to Y$ be the cover ramified along $S+T$ which we
considered in (\ref{diag1}). Choosing $U$ small enough,
we may assume that the connected component $U'\subset
\pi^{-1}(U)$ are polydisks with coordinates $\{w_1, \ldots , w_n\}$,
and that $\pi$ is defined by
$$
\pi(w_1,\ldots,w_n)=(z_1^{\mu_1},\ldots,z_n^{\mu_n}).
$$
Hence for $S'=\pi^*(S)_{red},$ and $T'=\pi^*(T)_{red},$
the restrictions to $U'$ are the zero sets of
$$
w_1\cdot \cdots \cdot w_{k} \mbox{ \ \ and \ \ }
w_{k+1}\cdot \cdots \cdot w_{k+k'},
$$
respectively.

Consider as above the Higgs bundle $\bigoplus {E'}^{p,q}$
obtained from the canonical extension of $\mathbb V_0'$ along
$T'$, and let $\{e'_1,\,e'_2, \, \ldots \}$ be a basis for
${E'}^{n-m,m}|_{U'}$.
\begin{claim}\label{7.7} For $U$ and $U'$ sufficiently small,
there exist some $\beta'\gg 1$ and a real number
$c>0$ with
$$
h'(e'_i(w),e'_j(w))\leq
c ((-\log|w_1|)\cdot(-\log|w_2|)\cdot
\cdots \cdot(-\log|w_{k}|))^{\beta'},
$$
for all $ w = (w_1,\ldots ,w_n) \in U'\setminus T'$.
\end{claim}
\begin{proof} By \cite{CKS}, 5.21, $U'_0=U'\setminus T'$ can be
decomposed into
$$
U'_0=\bigcup {U'}^I_{0,K},
$$
where the open subset $U_{0,K}^I$ depends on the index of the
filtration of the mixed Hodge structure (see Section 5.7 of
\cite{CKS}), and such that
$$
h'(e_i'(w),e_i'(w))\sim
(-\log|w_1|)^{l_1/2}\cdot(-\log|w_2|)^{(l_2-l_1)/2}\cdot \cdots
\cdot(-\log|w_{k}|)^{(l_{k}-l_{k-1})/2},
$$
for all $w\in {U'}^I_{0,K},$ where $(l_1,\, l_2,\ldots ,l_{k})$
is the multi index of the weight filtration of the mixed Hodge
structure. Since this index set is finite, there exist some
$\beta'\gg 1$ and some $c>0 $ such that
$$
h'(e_i'(w),e_i'(w))\leq c((-\log|w_1|)\cdot(-\log|w_2|)\cdot \cdots
\cdot(-\log|w_{k}|))^{\beta'},
$$
for all $w\in {U'}^I_{0,K}$ and for all $I$. Hence
$$
h'(e_i'(w),e_i'(w))\leq c((-\log|w_1|)\cdot(-\log|w_2|)\cdot\cdots
\cdot(-\log|w_{k}|))^{\beta'},
$$
for all $w\in U'\setminus T'$. By the Cauchy-Schwarz inequality we obtain
$$
h'(e'_i(w),e'_j(w))\leq
c((-\log|w_1|)\cdot(-\log|w_2|)\cdot\cdots\cdot(-\log|w_{k}|))^{\beta'},
$$
for all $w\in U'\setminus T'$.
\end{proof}

$Y$ is compact, hence there is a finite covering $\{U\}$ of $Y$
such that for all $U$ and each of the finitely many connected
components $U'$ of $\pi^{-1}(U)$, \ref{7.7} holds true.
We may even assume, that \ref{7.7} remains true, for the same
$\beta'$, for all point in a small neighborhood of the closure
$\bar{U}'$, not lying on $T'$.

We choose some $\alpha\gg 1$ such that for all the open sets $U'$
and for the constant $\beta'$ given by \ref{7.7}
one has
$$
\alpha \geq \beta' + 2m.
$$
In order to prove (\ref{equation1}) it is sufficient to show that
on each $U'$ there is some $\epsilon>0$ with
\begin{equation}\label{equation2}
||{\pi'}^* d\gamma(\partial_t)|_{{\gamma'}^{-1}(U')}
||^2_{{\pi'}^*\gamma^*\omega_{\alpha}}\geq
\epsilon {\pi'}^*\gamma^*(r^{-\alpha/m+2})||\tilde\tau((\partial_t)
^{m})|_{{\gamma'}^{-1}(U')}||^{\frac{2}{m}}_{{\pi'}^*\gamma^*(h\otimes
g^{-1})}.
\end{equation}

Let us return to the diagram (\ref{diag1}). As in the beginning
of this section, for each component $T_i$ of $T$ we consider
$\sL_i=\sO_Y(T_i)$ with the hermitian metric $g_i$, and
$\pi^*L_i$ with the pullback metric $\pi^*g_i.$ Let $s_i$ be a
section of $\sL_i$ with zero locus $T_i$, where we assume that
$s_i$ has been rescaled as needed in \ref{7.3} for the constant
$\alpha$, chosen above.

For the section $s'_i=\pi^*s_i$ define
$$
r'_i=-\log||s'_i||_{\pi*g_i}
$$
and $r'=r_1'\cdot\cdots\cdot r_\ell'.$ Obviously one has $r'_i=\pi^*r_i$ and
$r'=\pi^*r$.\\
\ \\
{\it Proof of the inequality} (\ref{equation2}).
Let $\{\phi_1,\phi_2, \ldots \} $ be an orthonormal basis for
$T_Y(-\log T)|_{U}$, with respect to $\omega_{\alpha}$. Then
$$
\{\phi_{i_1}\otimes\cdots\otimes\phi_{i_m}; \ i_1 \leq \cdots \leq
i_m\}
$$
is an orthonormal basis for $S^mT_Y(-\log T)|_{U}$ with respect
to $\omega_{\alpha}$ and
$$
\{\gamma^*(\phi_{i_1}\otimes\cdots\otimes\phi_{i_m}); \ i_1 \leq
\cdots \leq i_m\}
$$
is an orthonormal basis for $\gamma^*S^mT_Y(-\log T)|_{\gamma^{-1}U}$
with respect to $\gamma^*\omega_{\alpha}.$
Then, using the morphisms in (\ref{diag1}),
$$
\{{\gamma'}^{*}\pi^*(\phi_{i_1}\otimes\cdots \otimes \phi_{i_m});
\ i_1 \leq \cdots \leq i_m\}
$$
is an orthonormal basis for
${\gamma'}^{*}\pi^*S^mT_Y(-\log T)|_{{\gamma'}^{-1}U'}$ with respect
to ${\gamma'}^{*}\pi^*\omega_{\alpha}.$

For the map
\begin{equation}\label{map}
d\gamma^{m}: T^{m}_{\mathbb C}|_{\gamma^{-1}(U)}\to
\gamma^*(S^mT_{Y}(-\log T)|_{U}),
\end{equation}
write
$$
d\gamma^{m}((\partial_t)^{
m}|_{\gamma^{-1}(U)})=\sum
c_{i_1,\ldots,i_m}\gamma^*(\phi_{i_1}\otimes\cdots\otimes\phi_{i_m}).
$$
Then
$$
||d\gamma(\partial_t)|_{\gamma^{-1}(U)}||^2_{\gamma^*\omega_{\alpha}}=(\sum
|c_{i_1,\ldots,i_m}|^2)^{1/m}.
$$
Let
$$
{\pi'}^*d\gamma^{m}: {\pi'}^* T^{m}_{\mathbb C}|_{\gamma^{-1}(U)}\to
{\pi'}^*\gamma^*(S^mT_{Y}(-\log T)|_{U})
$$
be the pullback of the morphism (\ref{map}).
By the commutativity of (\ref{diag1}) one obtains
$$
{\pi'}^*d\gamma^{m}((\partial_t)^{
m})|_{{\gamma'}^{-1}{\pi}^{-1}(U)} =\sum
{\pi'}^*(c_{i_1,\ldots,i_m}){\gamma'}^{*}{\pi}^*(\phi_{i_1}\otimes
\cdots\otimes \phi_{i_m})
$$
and
$$
||{\pi'}^*d\gamma(\partial_t)|_{{\gamma'}^{-1}
\pi^{-1}(U)}||^2_{{\pi'}^*\gamma^*\omega_{\alpha}}
=(\sum {\pi'}^*|c_{i_1,\ldots,i_m}|^2)^{1/m}.
$$
Next we consider the second map
$$
\gamma^*\tau:\gamma^*(S^mT_Y(-\log T)|_U)\to
\gamma^*(E^{n-m,m}\otimes\sA^{-1}|_U)
$$
and its pullback
\begin{multline*}
{\pi'}^*\gamma^*\tau:{\gamma'}^{*}\pi^*(S^mT_Y(-\log
T)|_U)\to{\gamma'}^{*}\pi^*(E^{n-m,m}\otimes\sA^{-1}|_U)\\
\hookrightarrow{\gamma'}^{*}({E'}^{n-m,m}\otimes\pi^*\sA^{-1}|_{U'}).
\end{multline*}
For the connected component $U'$ of $\pi^{-1}(U)$ let
${a'}^{-1}$ be a local generator of $\pi^*\sA^{-1}|_{U'}.$ Then
$\{e'_1\otimes {a'}^{-1},\,e'_2\otimes {a'}^{-1},\ldots \}$ is a
basis of ${E'}^{n-m,m}\otimes\pi^*\sA^{-1}|_{{\gamma'}^{-1}(U)}$
and the morphism
$$
\pi^*\tau:\pi^*S^m T_Y(-\log T)|_{U'}\to
{E'}^{n-m,m}\otimes\pi^*\sA^{-1}|_{U'}
$$
is given by
$$
\pi^*\tau(\pi^*(\phi_{i_1}\otimes\cdots\otimes\phi_{i_m}))=
\sum b^j_{i_1,\ldots,i_m}e'_j\otimes {a'}^{-1}
$$
and one finds
$$
{\pi'}^*\gamma^*\tau d\gamma^{m}((\partial_t)^{m}
|_{{\gamma'}^{-1}(U')})=\sum {\pi'}^*(c_{i_1,\ldots,i_m})
{\gamma'}^{*}(b^j_{i_1,\ldots,i_m}){\gamma'}^{*}(e'_j\otimes {a'}^{-1}).
$$
Since the metric $\pi^*g^{-1}$ on $\pi^*\sA^{-1}$  is regular on
$U'$ the claim \ref{7.7} implies that
\begin{multline*}
|{\gamma'}^{*}(h'\otimes\pi^*g^{-1})({\gamma'}^{*}(e'_i\otimes
{a'}^{-1}), {\gamma'}^{*}(e'_j\otimes {a'}^{-1}))|\\
\leq c{\gamma'}^{*}((-\log|w_1|)\cdot
(-\log|w_2|)\cdot...\cdot(-\log|w_{k}|))^{\beta'}.
\end{multline*}
Here and later we allow ourselves to replace the constant $c$ by
some larger constant, whenever it is needed.

For the ramification order $\mu_i$ of $\pi$ over $T_i$,
and for some positive function $d_i$ on $U'$ one has
$$
|w_i|=d_i||{s'}_i^{\frac{1}{\mu_i}}|_{U'}||_{\pi^*g_i}.
$$
This description extends to the compactification $\bar{U}'$ of $U'$. Since
$\bar{U}'$ is compact, $d_i$ is bounded away from zero, and one finds
$$
|{\gamma'}^{*}(h'\otimes\pi^*g^{-1})({\gamma'}^{*}(e'_i\otimes
{a'}^{-1}), {\gamma'}^{*}(e'_j\otimes {a'}^{-1}))|
\leq c {\gamma'}^{*}{r'}^{\beta'}=c{\pi'}^*\gamma^*
r^{\beta'}.
$$
On the compact set $\bar U'$ all $b^j_{i_1,\ldots,i_m}$ are
bounded above. Hence, all ${\gamma'}^{*}(b^j_{i_1,\ldots,i_m})$
also are bounded above, and the Cauchy-Schwarz inequality implies
\begin{multline}\label{ineq}
||{\pi'}^*\tilde\tau
((\partial_t)^{m})|_{{\gamma'}^{-1}(U')}||^2_{{\pi'}^*\gamma^*(h\otimes
g^{-1})} \leq c{\pi'}^*\gamma^*r^{\beta'}\sum
{\pi'}^*|c_{i_1,\ldots,i_m}|^2\\
=c {\pi'}^*\gamma^*r^{\beta'}||{\pi'}^*d\gamma(\partial_t)
|_{{\gamma'}^{-1}(U')}||^{2m}_{{\pi'}^*\gamma^*\omega_\alpha}.
\end{multline}
Since we assumed $r \geq 1$ and $\alpha-2m \geq \beta'$, the
right hand side in (\ref{ineq}) is smaller than
$$
c {\pi'}^*\gamma^*r^{\alpha-2m}||{\pi'}^*d\gamma(\partial_t)
|_{{\gamma'}^{-1}(U')}||^{2m}_{{\pi'}^*\gamma^*\omega_\alpha},
$$
hence, we obtain the inequality
$$
||{\pi'}^*d\gamma(\partial_t)|_{{\gamma'}^{-1}(U')}||^2_{{\pi'}^*\gamma^*
\omega_{\alpha}}\geq \frac{1}{c}
{\pi'}^*\gamma^*(r^{-\alpha/m+2})||{\pi'}^*\tilde\tau((\partial_t)
^{m})|_{{\gamma'}^{-1}(U')}||^{\frac{2}{m}}_{{\pi'}^*
\gamma^*(h\otimes g^{-1})}
$$
as stated in (\ref{equation2}).
\end{proof}
{\it Proof of \ref{4.4}.} \ It remains to contradict the
existence of the ample sheaf $\sA$ and of the Higgs bundles
having the properties stated in \ref{6.2}. Those led to the
estimates in this section.

Recall the Ahlfors-Schwarz-Lemma, as stated in
1.1.1 in \cite{Siu} (see also lemma 3.2 in \cite{Dem}):
\begin{lemma}\label{ASL}
Let $c$ be a realvalued nonnegative function on $\C$ which locally is of the form $\varphi|f|^2$, where $\varphi$ is a local smooth positive function and $f$ is a local holomorphic function. Then there can not exists
any positive number $\rho$ such that
$$\partial_t \partial_{\bar{t}} \log c(t) \geq
\rho \cdot c(t)
$$
on $\C$ in the sense of currents.
\end{lemma}
Using the inequalities obtained in
\ref{7.10} and \ref{7.9} one has for suitable constants
$\epsilon$ and $\epsilon'$
\begin{multline*}
\frac{\sqrt{-1}}{2\pi} \partial_t \partial_{\bar{t}} \log c(t)
dt \wedge d \bar{t} = \frac{\sqrt{-1}}{2\pi} \partial \bar{\partial} \log
c(t) = - \Theta(T_{\mathbb C}, F)\\
\geq \epsilon \gamma^*\Theta(A_{Y\setminus T},
g_{\alpha}) \geq \epsilon\cdot \epsilon' \sqrt{-1} c(t)dt\wedge d\bar{t}
\end{multline*}
in the sense of currents. Hence
$$
\partial_t \partial_{\bar{t}} \log c(t) \geq
2\pi\cdot \epsilon \cdot \epsilon'\cdot c(t),
$$
contradicting Lemma \ref{ASL}.
\qed
%%%%%%%%%%%%%%%%%%%%%%%% References %%%%%%%%%%%%%%%%%
\bibliographystyle{plain}

\end{document}